Analytic Formulation and Simulation

# Mathematical analysis of 2D packing of circles on bounded and unbounded planes

*Harish Chandra Rajpoot*

*Research Scholar, Department of Mechanical Engineering, Indian Institute of Technology Bombay, Mumbai, Maharashtra, India-400076*

ABSTRACT

This paper encompasses the mathematical derivations of the analytic and generalized formula and recurrence relations to find out the radii of n umber of circles inscribed or packed in the plane region bounded by circular arcs (including sectors, semi and quarter circles) & the straight lines. The values of radii obtained using analytic formula and recurrence relations have been verified by comparing with those obtained using MATLAB codes. The methods used in this paper for packing circles are deterministic unlike heuristic strategies and optimization techniques. The analytic formulae derived for plane packing of tangent circles can be generalized and used for packing of spheres in 3D space and packing of circles on the spherical surface which is analogous to distribution of non-point charges. The packing density of identical circles, externally tangent to each other, the most densely packed on the regular hexagonal and the infinite planes have been formulated and analysed. This study paves the way for mathematically solving the problems of dense packing of circles in 2D containers, the packing of spheres in the voids (tetrahedral and octahedral) and finding the planar density on crystallographic plane.

*Keywords*: Tangent circles, Circle packing, Bounded plane, Packing density, Infinite plane

## 1. Introduction

In general, the circle packing is the arrangement of circles of equal or unequal radii on a surface (plane or curved, finite or infinite) without any overlap. The circle packing problems here are concerned with the packing of n number of external tangent circles in the plane region bounded by the circular arcs and the straight lines inside the square of known side, the sector of given radius and central angle and the regular hexagon of minimum size. The heuristic algorithms were used to pack the identical circles inside the plane container in form of triangle, square, rectangle, semicircle, and circle to maximize the radius for given size of container and to minimize the size of container for given radius [2]. The generic global optimization software packages were used to solve various circle packing problems and their industrial applications were reviewed [3]. The methods used here for circle packing are all deterministic and formulated analytically to compute the exact size of nth packed circle. The maximum density of the identical circles packed on the plane surface is obtained by arranging the centres of the circles in the hexagonal fashion which has been mathematically proved by theorem of Thue [4]. The most densely packed identical circles of a given radius can be enclosed by a regular hexagon of minimum side which can be determined in terms of number of circles and their radius.

Similarly, the three external tangent circles of given radii $a, b,$ and $c$ can be enclosed by a unique circle of minimum radius $R$ which has been analytically determined [1].

Conversely, if the radius of circumscribing circle $R$ and the radii say $a,$ and $b$ are known then the radius $c$ of third external tangent circle can be determined by solving the above equation for radius $c$ in terms of $a, b,$ and $R$.

Corresponding Author Phone: +91-9532487597
Email Address: hcrajpoot.iitb@gmail.com

## 1.1. Derivation of the radius of inscribed circle (one of three external tangent circles)

We know that the radius $(R)$ of circle which internally touches three external tangent circles with radii say $a, b$ & $c$ (as shown in figure-1) is given by the following generalized formula [1]

$$R = \frac{abc}{2\sqrt{abc(a+b+c)} - (ab+bc+ca)} \qquad (0 < a,b,c < R)$$

Where $a, b$ & $c$ $(a \geq b \geq c > 0)$ always satisfy the following inequality

$$c > \frac{ab}{(\sqrt{a} + \sqrt{b})^2}$$

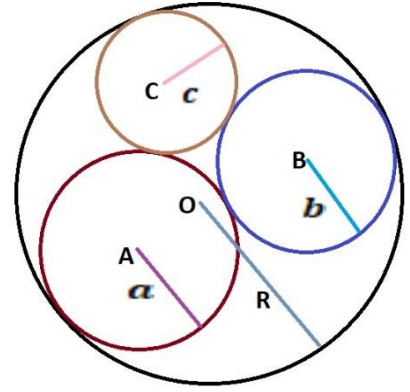

Now, suppose that the radius R of circumscribed circle with centre O and the radii of any two inscribed circles say $a$ & $b$ with centres A & B respectively are given or known then the radius $c$ of third inscribed circle with centre C (see figure-1) can be easily computed by solving above generalized equation for unknown radius $c$ as follows

$$R = \frac{abc}{2\sqrt{abc(a+b+c)} - (ab+bc+ca)}$$

**Figure-1:** Three externally touching circles with radii $a, b$ & $c$ are inscribed by the circle with center O and R

$$2R\sqrt{abc(a+b+c)} = abc + R(ab+bc+ca)$$

Taking the squares of both the positive sides of above equation,

$$\left(2R\sqrt{abc(a+b+c)}\right)^2 = \left(abc + R(ab+bc+ca)\right)^2$$

$$4R^2 ab(a+b)c + 4R^2 abc^2 = (aR+bR+ab)^2 c^2 + a^2 b^2 R^2 + 2abR(aR+bR+ab)c$$

$$((aR+bR+ab)^2 - 4abR^2)c^2 + \left(2abR(aR+bR+ab) - 4abR^2(a+b)\right)c + a^2 b^2 R^2 = 0$$

$$((aR+bR+ab)^2 - 4abR^2)c^2 - 2abR(aR+bR-ab)c + a^2 b^2 R^2 = 0$$

Solving above quadratic equation for $c$ using quadratic formula as follows

$$c = \frac{2abR(aR+bR-ab) \pm \sqrt{(2abR(aR+bR-ab))^2 - 4((aR+bR+ab)^2 - 4abR^2)a^2 b^2 R^2}}{2((aR+bR+ab)^2 - 4abR^2)}$$

$$= \frac{2abR\left(aR+bR-ab \pm \sqrt{(aR+bR-ab)^2 - ((aR+bR+ab)^2 - 4abR^2)}\right)}{2((aR+bR+ab)^2 - 4abR^2)}$$

$$= \frac{abR\left(aR+bR-ab \pm \sqrt{4abR^2 - 4a^2 bR - 4ab^2 R}\right)}{(aR+bR+ab)^2 - 4abR^2}$$

$$c = \frac{abR\left(aR+bR-ab \pm 2\sqrt{abR(R-a-b)}\right)}{(aR+bR+ab)^2 - 4abR^2}$$

In above expression it is worth noticing that $aR + bR - ab > 0$ & $(aR+bR+ab)^2 - 4abR^2 > 0$. The radius $c$ is a positive real value (i.e. can't be an imaginary number) therefore $(R - a - b) \geq 0$. Thus there are two possible positive real values of unknown radius $c$ which depends on the following two cases

Case-1: $R - a - b > 0$ or $R > a + b$    [discriminant is positive]

There are two positive distinct real values of unknown radius $c$ i.e. there are two different inscribed circles of radii $c_1$ & $c_2$ for given values of $a, b$ & $R$ (as shown in figure-2) given as follows

$$c_1, c_2 = \frac{abR\left(aR + bR - ab \pm 2\sqrt{abR(R - a - b)}\right)}{(aR + bR + ab)^2 - 4abR^2}$$

$\Rightarrow c_1 \neq c_2$

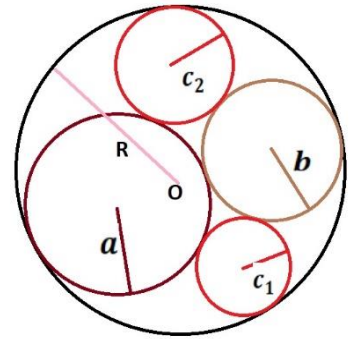

**Figure-2: Two different inscribed circles with radii $c_1$ & $c_2$ are possible if $R > a + b$. The circle with radius $c_1$ is inscribed/packed in the smaller region bounded by circles with radii $a, b$ and $R$.**

Case-2: $R - a - b = 0$ or $R = a + b$    [discriminant is zero]

There are two positive equal real values of unknown radius $c$ i.e. there are two identical inscribed circles of equal radii $c_1$ & $c_2$ for given values of $a, b$ & $R$. In this case, the centres A and B will lie on the diameter of outer/circumscribed circle with centre O (as shown in figure-3). The equal values of radius c are given as follows

$$c_1, c_2 = \frac{abR(aR + bR - ab)}{(aR + bR + ab)^2 - 4abR^2}$$

$\Rightarrow c_1 = c_2$

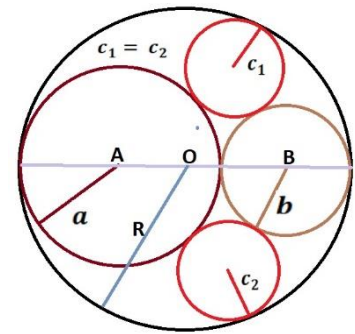

**Figure-3: Two identical inscribed circles with radii $c_1$ & $c_2$ are possible if $R = a + b$**

Case-3: $R - a - b < 0$ or $R < a + b$    [discriminant is negative]

In this case, there are no real values of $c_1$ & $c_2$ therefore this case is discarded.

Now, it's worth noticing that the value of radius $c$ of inscribed circle with centre C will be less than given radii $a, b$ & $R$ (see figure-3 or 4 in above two possible cases) only when we take the negative sign in the equation given by above case-1 (∵ -ve sign always gives decreasing value of radius $c$ of packed circle) i.e.

$$c_{\min} = \frac{abR\left(aR + bR - ab - 2\sqrt{abR(R - a - b)}\right)}{(aR + bR + ab)^2 - 4abR^2} \quad \forall \ c_{\min} < \min(a, b, R) \quad \ldots\ldots\ldots (1)$$

The above formula i.e. Eq(1) should be used for computing the minimum value of radius $c$ of (smallest) circle with the centre C which is inscribed/packed in the smaller region bounded by the tangent circles of given radii $a, b$ & $R$ (as shown in the above figure-2). The above eq(1) gives the decreasing value of $c$ which is applicable in all such cases of circle packing as discussed below.

## 2. Packing of circles in the plane bounded by square, semi and quarter circles

Consider a square ABCD of each side $x$ which has a semi-circle with diameter AB and a quarter circle with centre at the vertex B (see figure-4).

Now, consider the circles with centres $C_1, C, C_3, \ldots \ldots C_n$ and radii $r_1, r_2, r_3, \ldots \ldots r_n$ respectively snugly fitted or packed in the plane region bounded by the square ABCD, semi-circle and the quarter circle (as shown in the figure-4). Let's first find the radius $r_1$ of first inscribed (blue) circle then the radii of all other inscribed circles $r_2, r_3, \ldots \ldots r_n$ can be easily computed by using mathematical relation i.e. Eq(1) as derived above. We have to find out radius of nth inscribed circle (labelled with blue colour).

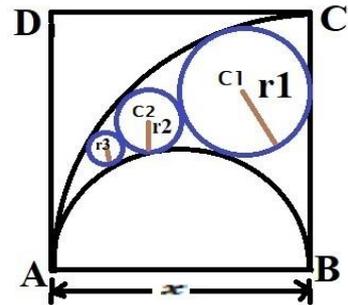

Figure-4: Infinite no. of small blue circles can be packed or inscribed in the plane region bounded by square, semi and quarter circles.

Now, join the centre $C_1$ of first blue inscribed circle to the centre E of semi-circle and to the vertex B of square ABCD. Drop the perpendiculars $C_1G$ and $C_1F$ from the centre $C_1$ to the sides BC and AB respectively (as shown in the figure-5 below)

In right $\Delta C_1GB$, We have, $C_1G = r_1$ and $C_1B = x - r_1$ (see figure-5). Now using Pythagorean Theorem as follows

$$(C_1B)^2 = (C_1G)^2 + (GB)^2$$

$$\Rightarrow \quad GB = \sqrt{(C_1B)^2 - (C_1G)^2}$$

$$GB = \sqrt{(x - r_1)^2 - (r_1)^2} = \sqrt{x^2 - 2xr_1}$$

In right $\Delta C_1FE$, We have (see figure-5)

$$C_1F = GB = \sqrt{x^2 - 2xr_1} \quad , \quad C_1E = \frac{x}{2} + r_1$$

$$FE = EB - FB = EB - C_1G$$
$$= \frac{x}{2} - r_1$$

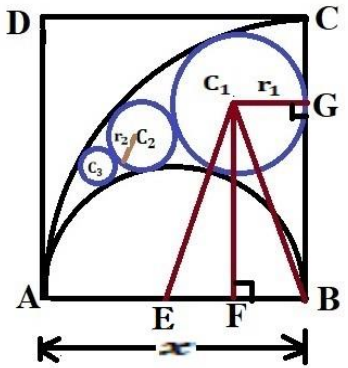

Figure-5: The blue circle of radius $r_1$ and center $C_1$ is the first largest circle inscribed in the region bounded by square ABCD, semi-circle with center E and quarter circle with center B.

Now using Pythagorean Theorem in right $\Delta C_1FE$ as follows

$$(C_1E)^2 = (C_1F)^2 + (FE)^2$$

$$\left(\frac{x}{2} + r_1\right)^2 = \left(\sqrt{x^2 - 2xr_1}\right)^2 + \left(\frac{x}{2} - r_1\right)^2$$

$$\frac{x^2}{4} + r_1^2 + xr_1 = x^2 - 2xr_1 + \frac{x^2}{4} + r_1^2 - xr_1$$

$$x^2 - 4xr_1 = 0 \quad \Rightarrow \quad r_1 = \frac{x}{4} \quad \forall \ x \neq 0$$

$$\boxed{r_1 = \frac{x}{4}} \qquad \ldots \ldots \ldots (2)$$

Above Eq(2) gives the radius of first largest inscribed circle (labelled with blue colour in the above figure-5) in terms of edge length $x$ of square ABCD.

**Coordinates of the centre ($C_1$):** The coordinates of the centre $C_1$ of first inscribed circle w.r.t. the vertex A (assumed to be the origin) of the square ABCD (see above figure-5) can be determined as follows

$$\text{x} - \text{coordinate} = \text{AF} = \text{AB} - \text{FB} = x - r_1 = x - \frac{x}{4} = \frac{3x}{4}$$

$$\text{y} - \text{coordinate} = C_1 F = GB = \sqrt{x^2 - 2xr_1} = \sqrt{x^2 - 2x \cdot \frac{x}{4}} = \frac{x}{\sqrt{2}}$$

Therefore, the coordinates of the centre $C_1$ of first inscribed circle w.r.t. the vertex A are $\left(\frac{3x}{4}, \frac{x}{\sqrt{2}}\right)$ (see the figure-5 above). Similarly, the coordinates of other inscribed circles can be computed after computing their radii.

**Radius of second inscribed circle $(r_2)$:** Now, draw the complete circles by extending the semi-circle with centre E and quarter-circle with centre B (as shown in the figure-6 below). It is worth noticing that three externally touching circles with centres $C_1$, $C_2$ and E touch the circumscribing circle (i.e. extended form of quarter circle in the above figure-5) with centre B. This case is similar to the case analysed in the above section 1.1 (see the above figure-1).

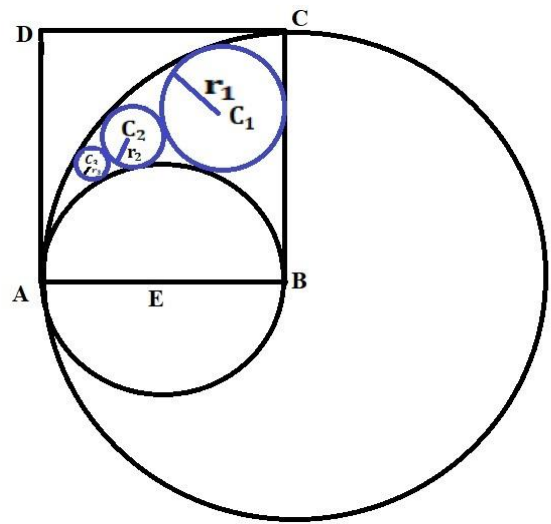

**Figure-6: The semi-circle and the quarter circle are extended to draw complete circles with centers E & B and radii $x/2$ & $x$ respectively.**

Now, the radius $r_2$ of the blue inscribed circle with centre $C_2$ touching the circles with centres $C_1$, E & B and the radii $r_1 = x/4$, $x/2$ & $x$ respectively, is found out by substituting $a = r_1 = x/4$, $b = x/2$ and $R = x$ in the above Eq(1) (i.e. minimum value $c_{\min}$) as follows

$$c_{\min} = \frac{abR\left(aR + bR - ab - 2\sqrt{abR(R - a - b)}\right)}{(aR + bR + ab)^2 - 4abR^2}$$

$$\Rightarrow r_2 = \frac{\frac{x}{4} \cdot \frac{x}{2} \cdot x \left(\frac{x}{4} \cdot x + \frac{x}{2} \cdot x - \frac{x}{4} \cdot \frac{x}{2} - 2\sqrt{\frac{x}{4} \cdot \frac{x}{2} \cdot x \left(x - \frac{x}{4} - \frac{x}{2}\right)}\right)}{\left(\frac{x}{4} \cdot x + \frac{x}{2} \cdot x + \frac{x}{4} \cdot \frac{x}{2}\right)^2 - 4 \cdot \frac{x}{4} \cdot \frac{x}{2} \cdot x^2}$$

$$= \frac{\frac{x^3}{8}\left(\frac{5x^2}{8} - \frac{2\sqrt{2}x^2}{8}\right)}{\left(\frac{7x^2}{8}\right)^2 - \frac{x^4}{2}} = \frac{\frac{x^5}{64}(5 - 2\sqrt{2})}{\frac{17x^4}{64}} = \frac{x}{17}(5 - 2\sqrt{2})$$

$$\therefore \quad r_2 = \frac{x}{17}\left(5 - 2\sqrt{2}\right) \approx 0.12773958089728293\, x \qquad \ldots\ldots\ldots (3)$$

Above Eq(3) gives the radius of second inscribed circle (labelled with blue colour in the above figure-6) in terms of edge length $x$ of square ABCD.

**Radius of third inscribed circle ($r_3$):** Similarly, consider three externally touching circles with centres $C_2$, $C_3$ and E which internally touch the largest circle with centre B (as shown in the figure-6 above). This case is again similar to the case analysed in above section 1.1. Thus, the value of radius $r_3$ of blue inscribed circle with centre $C_3$ depends on the radii $r_2, x/2$ & $x$ of circles with centres $C_2$, E & B respectively. Thus, the radius $r_3$ of circle with centre $C_3$ can be computed by substituting $a = r_2$, $b = x/2$ and $R = x$ in the above Eq(1) as follows

$$c_{\min} = \frac{abR\left(aR + bR - ab - 2\sqrt{abR(R-a-b)}\right)}{(aR+bR+ab)^2 - 4abR^2}$$

$\Rightarrow r_3$

$$= \frac{\frac{x}{17}(5-2\sqrt{2}) \cdot \frac{x}{2} \cdot x \left(\frac{x}{17}(5-2\sqrt{2}) \cdot x + \frac{x}{2} \cdot x - \frac{x}{17}(5-2\sqrt{2}) \cdot \frac{x}{2} - 2\sqrt{\frac{x}{17}(5-2\sqrt{2}) \cdot \frac{x}{2} \cdot x \left(x - \frac{x}{17}(5-2\sqrt{2}) - \frac{x}{2}\right)}\right)}{\left(\frac{x}{17}(5-2\sqrt{2}) \cdot x + \frac{x}{2} \cdot x + \frac{x}{17}(5-2\sqrt{2}) \cdot \frac{x}{2}\right)^2 - 4 \cdot \frac{x}{17}(5-2\sqrt{2}) \cdot \frac{x}{2} \cdot x^2}$$

$$= \frac{x\left(33 - 20\sqrt{2} + 85 - 34\sqrt{2} - 2(5-2\sqrt{2})\sqrt{85 - 34\sqrt{2} - 66 + 40\sqrt{2}}\right)}{289 + 297 - 180\sqrt{2} - 170 + 68\sqrt{2}}$$

$$= \frac{x\left(118 - 54\sqrt{2} - 2(5-2\sqrt{2})(3\sqrt{2}+1)\right)}{16(26 - 7\sqrt{2})}$$

$$= \frac{x(33 - 20\sqrt{2})}{4(26 - 7\sqrt{2})} = \frac{x(33-20\sqrt{2})(26+7\sqrt{2})}{4(26-7\sqrt{2})(26+7\sqrt{2})} = \frac{x(578 - 289\sqrt{2})}{4(578)} = \frac{x(2-\sqrt{2})}{8}$$

$$\therefore \quad r_3 = \frac{x}{8}(2-\sqrt{2}) \approx 0.0732233047033631 x \qquad \dots \dots \dots (4)$$

Above Eq(4) gives the radius of third inscribed circle (labelled with blue colour in the above figure-6) in terms of edge length $x$ of square ABCD.

**Radius of forth inscribed circle ($r_4$):** Similarly, the radius $r_4$ can be computed by substituting $a = r_3$, $b = x/2$ and $R = x$ in the above Eq(1) as follows

$$c_{\min} = \frac{abR\left(aR + bR - ab - 2\sqrt{abR(R-a-b)}\right)}{(aR+bR+ab)^2 - 4abR^2}$$

$\Rightarrow r_4$

$$= \frac{\frac{x}{8}(2-\sqrt{2}) \cdot \frac{x}{2} \cdot x \left(\frac{x}{8}(2-\sqrt{2}) \cdot x + \frac{x}{2} \cdot x - \frac{x}{8}(2-\sqrt{2}) \cdot \frac{x}{2} - 2\sqrt{\frac{x}{8}(2-\sqrt{2}) \cdot \frac{x}{2} \cdot x \left(x - \frac{x}{8}(2-\sqrt{2}) - \frac{x}{2}\right)}\right)}{\left(\frac{x}{8}(2-\sqrt{2}) \cdot x + \frac{x}{2} \cdot x + \frac{x}{8}(2-\sqrt{2}) \cdot \frac{x}{2}\right)^2 - 4 \cdot \frac{x}{8}(2-\sqrt{2}) \cdot \frac{x}{2} \cdot x^2}$$

$$= \frac{x\left(6 - 4\sqrt{2} + 16 - 8\sqrt{2} - 2(2-\sqrt{2})\sqrt{16 - 8\sqrt{2} - 12 + 8\sqrt{2}}\right)}{64 + 54 - 36\sqrt{2} - 32 + 16\sqrt{2}}$$

$$= \frac{x\left(22 - 12\sqrt{2} - 2(2-\sqrt{2})(2)\right)}{2(43 - 10\sqrt{2})} = \frac{x(7-4\sqrt{2})}{43 - 10\sqrt{2}} = \frac{x(7-4\sqrt{2})(43+10\sqrt{2})}{(43-10\sqrt{2})(43+10\sqrt{2})}$$

$$= \frac{x(301 - 172\sqrt{2} + 70\sqrt{2} - 80)}{1849 - 200} = \frac{x(221 - 102\sqrt{2})}{1649}$$

$$\therefore \; r_4 = \frac{x(221 - 102\sqrt{2})}{1649} \approx 0.046543490987231234x \qquad \ldots \ldots \ldots (5)$$

Above Eq(5) gives the radius of fourth inscribed circle (see the above figure-6) in terms of edge length $x$ of square ABCD.

**Radius of fifth inscribed circle ($r_5$):** Similarly, the radius $r_5$ can be computed by substituting $a = r_4$, $b = x/2$ and $R = x$ in the above Eq(1) as follows

$$c_{\min} = \frac{abR\left(aR + bR - ab - 2\sqrt{abR(R - a - b)}\right)}{(aR + bR + ab)^2 - 4abR^2}$$

$\Rightarrow r_5$

$$= \frac{\frac{x(221 - 102\sqrt{2})}{1649} \cdot \frac{x}{2} \cdot x \left(\frac{x(221 - 102\sqrt{2})}{1649} \cdot x + \frac{x}{2} \cdot x - \frac{x(221 - 102\sqrt{2})}{1649} \cdot \frac{x}{2} - 2\sqrt{\frac{x(221 - 102\sqrt{2})}{1649} \cdot \frac{x}{2} \cdot x \left(x - \frac{x(221 - 102\sqrt{2})}{1649} - \frac{x}{2}\right)}\right)}{\left(\frac{x(221 - 102\sqrt{2})}{1649} \cdot x + \frac{x}{2} \cdot x + \frac{x(221 - 102\sqrt{2})}{1649} \cdot \frac{x}{2}\right)^2 - 4 \cdot \frac{x(221 - 102\sqrt{2})}{1649} \cdot \frac{x}{2} \cdot x^2}$$

$$= \frac{x\left(69649 - 45084\sqrt{2} + 364429 - 168198\sqrt{2} - 2(221 - 102\sqrt{2})\sqrt{364429 - 168198\sqrt{2} - 139298 + 90168\sqrt{2}}\right)}{2719201 + 626841 - 405756\sqrt{2} - 728858 + 336396\sqrt{2}}$$

$$= \frac{x\left(434078 - 213282\sqrt{2} - 2(221 - 102\sqrt{2})(459 - 85\sqrt{2})\right)}{2617184 - 69360\sqrt{2}} = \frac{x(196520 - 82076\sqrt{2})}{4624(566 - 15\sqrt{2})}$$

$$= \frac{x(170 - 71\sqrt{2})}{4(566 - 15\sqrt{2})} = \frac{x(170 - 71\sqrt{2})(566 + 15\sqrt{2})}{4(566 - 15\sqrt{2})(566 + 15\sqrt{2})} = \frac{x(47045 - 18818\sqrt{2})}{639812}$$

$$\therefore \; r_5 = \frac{x(47045 - 18818\sqrt{2})}{639812} \approx 0.03193489522432073x \qquad \ldots \ldots \ldots (6)$$

Above Eq(6) gives the radius of fifth inscribed circle (see the above figure-6) in terms of edge length $x$ of square ABCD.

### 3.1 Recurrence Relation

Similarly, the radius $r_6$ of sixth inscribed circle can be computed by substituting $a = r_5$, $b = x/2$ and $R = x$ in the above Eq(1) and so on. Every time the case becomes similar to the case of three externally touching circles which internally touch the fourth circumscribing circle which is analysed in the above section 1.1. This way we can compute the radius of any nth blue circle inscribed in the bounded plane region (see above figure-6). It's worth noticing that the value of radii $b = x/2$ and $R = x$ are constant (while radius $a$ is variable) in each case and the radius $r_n$ of nth inscribed circle depends on the value of radius $r_{n-1}$ of (n-1)th inscribed circle.

It can be inferred from the above derivations that the radius $r_{n+1}$ of (n+1)th circle inscribed in the bounded region can be computed by substituting $a = r_n$, $b = x/2$ and $R = x$ in the above Eq(1) as follows

$$c_{\min} = \frac{abR\left(aR + bR - ab - 2\sqrt{abR(R - a - b)}\right)}{(aR + bR + ab)^2 - 4abR^2}$$

$$\Rightarrow r_{n+1} = \frac{r_n \cdot \frac{x}{2} \cdot x \left( r_n \cdot x + \frac{x}{2} \cdot x - r_n \cdot \frac{x}{2} - 2\sqrt{r_n \cdot \frac{x}{2} \cdot x \left( x - r_n - \frac{x}{2} \right)} \right)}{\left( r_n \cdot x + \frac{x}{2} \cdot x + r_n \cdot \frac{x}{2} \right)^2 - 4 \cdot r_n \cdot \frac{x}{2} \cdot x^2}$$

$$= \frac{\frac{1}{2} x^2 r_n \left( \frac{1}{2} x r_n + \frac{x^2}{2} - 2\sqrt{\frac{1}{2} x^2 r_n \left( \frac{x}{2} - r_n \right)} \right)}{\left( \frac{3}{2} x r_n + \frac{x^2}{2} \right)^2 - 2 x^3 r_n}$$

$$= \frac{\frac{1}{4} x^3 r_n \left( r_n + x - 2\sqrt{r_n (x - 2 r_n)} \right)}{\frac{x^2}{4} (3 r_n + x)^2 - 2 x^3 r_n}$$

$$\therefore \; r_{n+1} = x \left( \frac{r_n^2 + x r_n - 2 r_n \sqrt{x r_n - 2 r_n^2}}{9 r_n^2 - 2 x r_n + x^2} \right) \quad \ldots\ldots\ldots\ldots (7)$$

Where, $r_1 = \frac{x}{4} \quad \forall \; x > 0, \; x \in R, \; n \in N \quad \Rightarrow \quad r_{n+1} = f(r_n)$

The above Eq(7) is the **governing equation** i.e. **recurrence relation** for 2D packing of circles to analytically compute the radius of (n+1)th inscribed circle by substituting the value of radius of nth inscribed circle (i.e. blue inscribed circles as shown in the above figure-6). It is worth noticing from above Eq(7) that next term $r_{n+1}$ is a function of the previous term $r_n$ for given side $x$ of the square.

### 3.2 MATLAB code-1

It is usually very difficult to compute the radius of nth (blue) inscribed circle when $n$ is very large therefore using the above recurrence relation i.e. Eq(7), a for-loop can easily be used to develop a simple MATLAB code, to compute the radii of $n$ number of inscribed circles when the edge length $x$ of square is given, as follows

```
clc;
clear all;
close all;
x=input('Enter the edge or side of square (in required unit) = ');
n=input('Enter the number of circles to be packed, n = ');
disp('Column matrix of radii of n no. of packed/inscribed circles, r(i)(units)=');
r(1)=x/4;
disp(r(1));
for i=1:1:n-1
    r(1)=x/4;
%HCR's recurrence relation to compute radii of inscribed circles
r(i+1)=double(x*((r(i))^2+x*r(i)-2*r(i)*sqrt(x*r(i)-2*(r(i))^2))/(9*(r(i))^2-
2*x*r(i)+x^2));
A=[r(i+1)];
disp(A);
end
```

Above MATLAB code of recurrence relation i.e. Eq(7) requires two input values 1) side $x$ of square and 2) number $n$ of circles to be inscribed/packed to compute the radii of all n number of circles packed or inscribed in the plane region bounded by given square, semi-circle and quarter circle as shown in the figure-5 above.

## 3.3 Results

After entering the values of side of unit square i.e. $x = 1$ and number of circles to packed/inscribed i.e. $n = 100$, the following are the radii of all 100 packed/inscribed as tabulated below.

**Inputs**

```
Enter the edge or side of square (in required unit) =
1
Enter the number of circles to be packed, n =
100
```

**Outputs**

```
Column matrix of radii of n no. of packed/inscribed circles, r(i)(units)=
```

**Table-1: Ratio of radius of packed/inscribed circle to the side of square (n=100 circles)**

| | | | |
|---|---|---|---|
| 0.2500 | 0.0014 | 3.7801e-04 | 1.7120e-04 |
| 0.1277 | 0.0013 | 3.6374e-04 | 1.6681e-04 |
| 0.0732 | 0.0012 | 3.5025e-04 | 1.6258e-04 |
| 0.0465 | 0.0012 | 3.3751e-04 | 1.5851e-04 |
| 0.0319 | 0.0011 | 3.2544e-04 | 1.5460e-04 |
| 0.0232 | 0.0010 | 3.1401e-04 | 1.5082e-04 |
| 0.0176 | 9.4996e-04 | 3.0318e-04 | 1.4719e-04 |
| 0.0137 | 8.9405e-04 | 2.9289e-04 | 1.4368e-04 |
| 0.0110 | 8.4293e-04 | 2.8312e-04 | 1.4030e-04 |
| 0.0091 | 7.9607e-04 | 2.7383e-04 | 1.3703e-04 |
| 0.0076 | 7.5301e-04 | 2.6499e-04 | 1.3388e-04 |
| 0.0064 | 7.1336e-04 | 2.5657e-04 | 1.3083e-04 |
| 0.0055 | 6.7675e-04 | 2.4855e-04 | 1.2789e-04 |
| 0.0048 | 6.4289e-04 | 2.4089e-04 | 1.2505e-04 |
| 0.0042 | 6.1151e-04 | 2.3359e-04 | 1.2230e-04 |
| 0.0037 | 5.8236e-04 | 2.2661e-04 | 1.1964e-04 |
| 0.0033 | 5.5526e-04 | 2.1994e-04 | 1.1706e-04 |
| 0.0029 | 5.3000e-04 | 2.1356e-04 | 1.1457e-04 |
| 0.0026 | 5.0643e-04 | 2.0745e-04 | 1.1216e-04 |
| 0.0024 | 4.8439e-04 | 2.0161e-04 | 1.0982e-04 |
| 0.0022 | 4.6376e-04 | 1.9600e-04 | 1.0755e-04 |
| 0.0020 | 4.4442e-04 | 1.9063e-04 | 1.0536e-04 |
| 0.0018 | 4.2627e-04 | 1.8547e-04 | 1.0323e-04 |
| 0.0017 | 4.0920e-04 | 1.8052e-04 | 1.0116e-04 |
| 0.0015 | 3.9314e-04 | 1.7577e-04 | 9.9157e-05 |

It's worth noticing that each value of above Table-1 shows the ratio of radius of a packed/inscribed circle to the side of square (i.e. taken as unity). Therefore, the exact value of radius of any of packed/inscribed circle of any order (i.e. up to $n = 100$) can be obtained by directly multiplying given side of square $x$ with the corresponding value of ratio as obtained in the above Table-1 (each column has 25 values of ratio) for 100 packed circles.

**Note:** The output values of radii of inscribed circles are obtained in form of a column from the above MALAB code-1 which have been rearranged in sequence in form of above Table-1 with four columns.

## 4. Packing of circles in the plane bounded by two semi circles and a quarter circle

Consider a square ABCD of each side $x$ which has two semi-circles with diameter AB & AD and a quarter circle with centre at the vertex B (see figure-7).

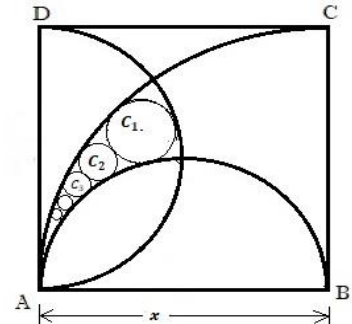

Now, consider the circles with centres $C_1, C_2, C_3, \ldots \ldots C_n$ and radii $r_1, r_2, r_3, \ldots \ldots r_n$ respectively snugly fitted or packed in the plane region bounded by two semicircles and a quarter circle inside the square ABCD (as shown in the figure-7). Let's first find the radius $r_1$ of first packed/inscribed circle then the radii of all other inscribed circles $r_2, r_3, \ldots \ldots r_n$ can be easily computed by using mathematical relation i.e. Eq(1) similar to the above section 3. We have to find out the radius of $n$th inscribed circle using recurrence relation i.e. Eq(7) as derived in the above section 3.

**Figure-7:** Infinite no. of small circles can be packed/inscribed in the region bounded by two semicircles and a quarter circle inside the square ABCD

Now, join the centre $C_1$ of first packed/inscribed circle to the centre E & G of semi-circles and to the vertex B of square ABCD. Drop the perpendiculars $C_1F$ and $C_1H$ from the centre $C_1$ to the sides AB and AD respectively (as shown in the figure-8 below). Let $C_1F = y$.

In right $\Delta C_1FB$, We have, $C_1F = y$ and $C_1B = x - r_1$ (see figure-8 below). Now using Pythagorean Theorem as follows

$$(C_1B)^2 = (C_1F)^2 + (FB)^2$$

$$\Rightarrow FB = \sqrt{(C_1B)^2 - (C_1F)^2}$$

$$FB = \sqrt{(x - r_1)^2 - y^2}$$

In right $\Delta C_1FE$, We have, $C_1F = y$, $C_1E = \frac{x}{2} + r_1$ (see figure-8)

Now using Pythagorean Theorem as follows

$$(C_1E)^2 = (C_1F)^2 + (FE)^2 \Rightarrow FE = \sqrt{(C_1E)^2 - (C_1F)^2}$$

$$FE = \sqrt{\left(\frac{x}{2} + r_1\right)^2 - y^2}$$

From figure-8, we have

$$FB = FE + EB \Leftrightarrow (FB)^2 = (FE + EB)^2$$

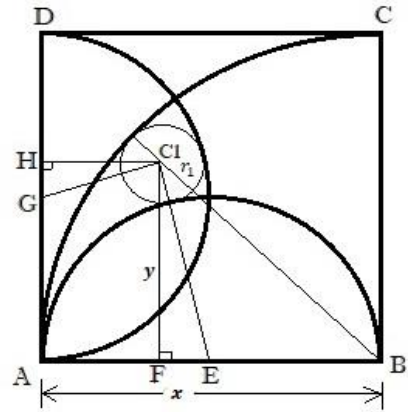

**Figure-8:** A circle of radius $r_1$ and center $C_1$ is the first largest circle inscribed in the region bounded by two semi-circles and a quarter circle inside the square ABCD.

$$\Rightarrow \left(\sqrt{(x-r_1)^2 - y^2}\right)^2 = \left(\sqrt{\left(\frac{x}{2} + r_1\right)^2 - y^2} + \frac{x}{2}\right)^2$$

$$x^2 + r_1{}^2 - 2xr_1 - y^2 = \frac{x^2}{4} + r_1{}^2 + xr_1 - y^2 + \frac{x^2}{4} + x\sqrt{\left(\frac{x}{2} + r_1\right)^2 - y^2}$$

$$\frac{x^2}{2} - 3xr_1 - x\sqrt{\left(\frac{x}{2} + r_1\right)^2 - y^2} = 0$$

$$x\left(\frac{x}{2} - 3r_1 - \sqrt{\left(\frac{x}{2} + r_1\right)^2 - y^2}\right) = 0$$

$$\frac{x}{2} - 3r_1 = \sqrt{\left(\frac{x}{2} + r_1\right)^2 - y^2} \qquad (\because\ x \ne 0)$$

$$\Rightarrow \left(\frac{x}{2} - 3r_1\right)^2 = \left(\sqrt{\left(\frac{x}{2} + r_1\right)^2 - y^2}\right)^2 \qquad (\because \text{LHS} > 0,\ \text{RHS} > 0)$$

$$\frac{x^2}{4} + 9r_1{}^2 - 3xr_1 = \frac{x^2}{4} + r_1{}^2 + xr_1 - y^2$$

$$y^2 = 4xr_1 - 8r_1{}^2 \qquad \ldots\ldots\ldots\ldots (4.1)$$

In right $\Delta C_1HG$, We have,

$$C_1G = \frac{x}{2} - r_1,\ \ C_1H = AF = AE - FE = \frac{x}{2} - \sqrt{\left(\frac{x}{2} + r_1\right)^2 - y^2}\ \ \text{and}\ \ HG = HA - GA = C_1F - GA = y - \frac{x}{2}$$

Now using Pythagorean Theorem (see figure-8 above) as follows

$$(C_1G)^2 = (C_1H)^2 + (HG)^2$$

$$\Rightarrow \left(\frac{x}{2} - r_1\right)^2 = \left(\frac{x}{2} - \sqrt{\left(\frac{x}{2} + r_1\right)^2 - y^2}\right)^2 + \left(y - \frac{x}{2}\right)^2$$

$$\frac{x^2}{4} + r_1^2 - xr_1 = \frac{x^2}{4} + \frac{x^2}{4} + r_1^2 + xr_1 - y^2 - x\sqrt{\left(\frac{x}{2} + r_1\right)^2 - y^2} + y^2 + \frac{x^2}{4} - xy$$

$$x\left(\frac{x}{2} + 2r_1 - y - \sqrt{\left(\frac{x}{2} + r_1\right)^2 - y^2}\right) = 0$$

$$\frac{x}{2} + 2r_1 - y = \sqrt{\left(\frac{x}{2} + r_1\right)^2 - y^2} \qquad (\because x \neq 0)$$

$$\Rightarrow \left(\frac{x}{2} + 2r_1 - y\right)^2 = \left(\sqrt{\left(\frac{x}{2} + r_1\right)^2 - y^2}\right)^2 \qquad (\because \text{LHS} > 0, \; \text{RHS} > 0)$$

$$\frac{x^2}{4} + 4r_1^2 + y^2 + 2xr_1 - 4r_1 y - xy = \frac{x^2}{4} + r_1^2 + xr_1 - y^2$$

$$2y^2 - (x + 4r_1)y + xr_1 + 3r_1^2 = 0$$

$$2(4xr_1 - 8r_1^2) - (x + 4r_1)y + xr_1 + 3r_1^2 = 0 \qquad \text{(Setting value of } y^2 \text{ from Eq(4.1))}$$

$$y = \frac{9xr_1 - 13r_1^2}{x + 4r_1} \qquad \ldots\ldots\ldots (4.2)$$

Now, substituting the value of $y$ from Eq(4.2) into above Eq(4.1) as follows

$$\left(\frac{9xr_1 - 13r_1^2}{x + 4r_1}\right)^2 = 4xr_1 - 8r_1^2$$

$$\frac{r_1^2(9x - 13r_1)^2 - 4r_1(x - 2r_1)(x + 4r_1)^2}{(x + 4r_1)^2} = 0$$

$$r_1\{r_1(81x^2 + 169r_1^2 - 234xr_1) - (4x - 8r_1)(x^2 + 16r_1^2 + 8xr_1)\} = 0 \qquad (\because x + 4r_1 \neq 0)$$

$$81x^2 r_1 + 169r_1^3 - 234xr_1^2 - 4x^3 - 64xr_1^2 - 32x^2 r_1 + 8x^2 r_1 + 128r_1^3 + 64xr_1^2 = 0 \qquad (\because r_1 \neq 0)$$

$$297r_1^3 - 234xr_1^2 + 57x^2 r_1 - 4x^3 = 0$$

$$297\left(\frac{r_1}{x}\right)^3 - 234\left(\frac{r_1}{x}\right)^2 + 57\left(\frac{r_1}{x}\right) - 4 = 0 \qquad \text{(Dividing by } x^3\text{)} \qquad \ldots (4.3)$$

The above Eq(4.3) is a cubic equation in terms of $\frac{r_1}{x}$ which can be easily solved by numerical techniques which give the following solutions

$$\frac{r_1}{x} = \frac{1}{3}, \frac{1}{3}, \frac{4}{33} \quad \Rightarrow \quad r_1 = \frac{x}{3}, \frac{x}{3}, \frac{4x}{33}$$

It is to be noted that the region bounded by two semicircles and a quarter circle inside the square ABCD (see the above figure-7) is less than the region bounded by semicircle, quarter circle and the square ABCD (see the above figure-4). This implies that the radius $r_1$ of first inscribed circle in this case (section 4) must be less than the radius $r_1 = \frac{x}{4}$ obtained in the above section 3. Therefore, We select the real roots of above cubic equation (4.3) such that $\frac{r_1}{x} < \frac{1}{4}$ or $r_1 < \frac{x}{4}$. After applying this condition, we obtain a unique solution as follows

$$\frac{r_1}{x} = \frac{4}{33} \quad \Rightarrow \quad r_1 = \frac{4x}{33} \qquad \left(\because r_1 < \frac{x}{4} \text{ or } r_1 \neq \frac{x}{3}\right)$$

$$\therefore r_1 = \frac{4x}{33} \qquad \ldots\ldots\ldots\ldots\ldots(8)$$

The above Eq(8) gives the radius of first largest packed/inscribed circle (as shown in the above figure-7) in terms of edge length $x$ of square ABCD.

**Coordinates of the centre ($C_1$):** The coordinates of the centre $C_1$ of first inscribed circle w.r.t. the vertex A (assumed to be the origin) of the square ABCD (see above figure-8) can be determined as follows

$$\text{y} - \text{coordinate} = C_1 F = y = \frac{9xr_1 - 13r_1{}^2}{x + 4r_1} = \frac{9x \cdot \frac{4x}{33} - 13\left(\frac{4x}{33}\right)^2}{x + 4 \cdot \frac{4x}{33}} = \frac{980x}{1617} = \frac{20x}{33}$$

$$\text{x} - \text{coordinate} = C_1 H = AF = AE - FE = \frac{x}{2} - \sqrt{\left(\frac{x}{2} + r_1\right)^2 - y^2} = \frac{x}{2} - \sqrt{\left(\frac{x}{2} + \frac{4x}{33}\right)^2 - \left(\frac{20x}{33}\right)^2}$$

$$= \frac{4x}{11}$$

Therefore, the coordinates of the centre $C_1$ of first inscribed circle w.r.t. the vertex A are $\left(\frac{4x}{11}, \frac{20x}{33}\right)$ (see the figure-8 above). Similarly, the coordinates of other inscribed circles can be computed after computing their radii.

**Radius of second inscribed circle ($r_2$):** Now, draw the complete circles by extending the semi-circle with centre E and quarter-circle with centre B (as shown in the figure-9 below). It is worth noticing that three externally touching circles with centres $C_1$, $C_2$ and E touch the circumscribing circle (i.e. extended form of quarter circle in the above figure-7) with centre B. This case is similar to the case of three externally touching circles which touch the forth circumscribing circle analysed in the above section 1.1 (see the above figure-1).

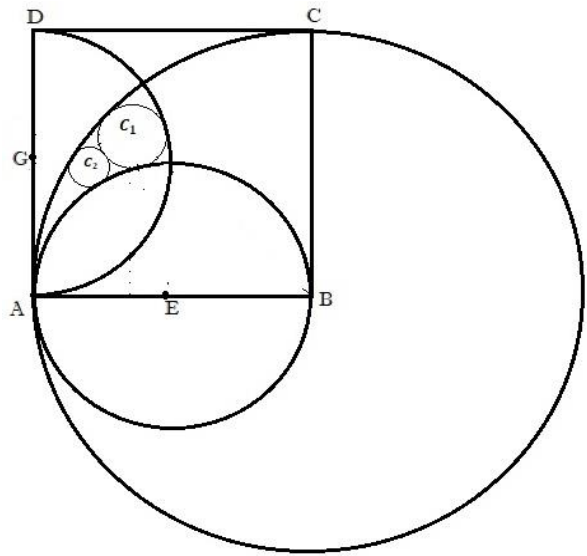

**Figure-9: The semi-circle and the quarter circle are extended to draw complete circles with centers E & B and radii $x/2$ & $x$ respectively.**

Now, the radius $r_2$ of the second inscribed circle with centre $C_2$ touching the circles with centres $C_1$, E & B and the radii $r_1 = 4x/33$, $x/2$ & $x$ respectively, is found out by substituting $a = r_1 = 4x/33$, $b = x/2$ and $R = x$ in the above Eq(1) (i.e. minimum value $c_{\min}$) as follows

$$c_{\min} = \frac{abR\left(aR + bR - ab - 2\sqrt{abR(R-a-b)}\right)}{(aR + bR + ab)^2 - 4abR^2}$$

$$\Rightarrow r_2 = \frac{\frac{4x}{33} \cdot \frac{x}{2} \cdot x \left(\frac{4x}{33} \cdot x + \frac{x}{2} \cdot x - \frac{4x}{33} \cdot \frac{x}{2} - 2\sqrt{\frac{4x}{33} \cdot \frac{x}{2} \cdot x \left(x - \frac{4x}{33} - \frac{x}{2}\right)}\right)}{\left(\frac{4x}{33} \cdot x + \frac{x}{2} \cdot x + \frac{4x}{33} \cdot \frac{x}{2}\right)^2 - 4 \cdot \frac{4x}{33} \cdot \frac{x}{2} \cdot x^2}$$

$$= \frac{\frac{2x^3}{33}\left(\frac{37x^2}{66} - \frac{10x^2}{33}\right)}{\left(\frac{15x^2}{22}\right)^2 - \frac{8x^4}{33}} = \frac{\frac{68x^5}{66^2}}{\frac{969x^4}{66^2}} = \frac{4x}{57}$$

$$\therefore \quad r_2 = \frac{4x}{57} \approx 0.07017543859649122\, x \qquad \ldots\ldots\ldots (9)$$

The above Eq(9) gives the radius of second inscribed circle (as shown in the above figure-7) in terms of edge length $x$ of square ABCD.

**Radius of third inscribed circle $(r_3)$:** Similarly, consider three externally touching circles with centres $C_2$, $C_3$ and E which internally touch the largest circle with centre B (as shown in the figure-7 above). This case is again similar to the case analysed in the above section 1.1. Thus, the radius $r_3$ of circle with centre $C_3$ can be computed by substituting $a = r_2$, $b = x/2$ and $R = x$ in the above Eq(1) as follows

$$c_{\min} = \frac{abR\left(aR + bR - ab - 2\sqrt{abR(R-a-b)}\right)}{(aR + bR + ab)^2 - 4abR^2}$$

$$\Rightarrow r_3 = \frac{\frac{4x}{57} \cdot \frac{x}{2} \cdot x \left(\frac{4x}{57} \cdot x + \frac{x}{2} \cdot x - \frac{4x}{57} \cdot \frac{x}{2} - 2\sqrt{\frac{4x}{57} \cdot \frac{x}{2} \cdot x \left(x - \frac{4x}{57} - \frac{x}{2}\right)}\right)}{\left(\frac{4x}{57} \cdot x + \frac{x}{2} \cdot x + \frac{4x}{57} \cdot \frac{x}{2}\right)^2 - 4 \cdot \frac{4x}{57} \cdot \frac{x}{2} \cdot x^2}$$

$$= \frac{\frac{2x^3}{57}\left(\frac{61x^2}{114} - \frac{14x^2}{57}\right)}{\left(\frac{69x^2}{114}\right)^2 - \frac{8x^4}{57}} = \frac{\frac{66x^5}{114^2}}{\frac{2937x^4}{114^2}} = \frac{4x}{89}$$

$$\therefore \quad r_3 = \frac{4x}{89} \approx 0.0449438202247191\, x \qquad \ldots\ldots\ldots (10)$$

The above Eq(10) gives the radius of third inscribed circle (as shown in the above figure-7) in terms of edge length $x$ of square ABCD.

**Radius of forth inscribed circle $(r_4)$:** Similarly, the radius $r_4$ can be computed by substituting $a = r_3$, $b = x/2$ and $R = x$ in the above Eq(1) as follows

$$c_{min} = \frac{abR\left(aR + bR - ab - 2\sqrt{abR(R - a - b)}\right)}{(aR + bR + ab)^2 - 4abR^2}$$

$$\Rightarrow r_4 = \frac{\frac{4x}{89} \cdot \frac{x}{2} \cdot x \left(\frac{4x}{89} \cdot x + \frac{x}{2} \cdot x - \frac{4x}{89} \cdot \frac{x}{2} - 2\sqrt{\frac{4x}{89} \cdot \frac{x}{2} \cdot x \left(x - \frac{4x}{89} - \frac{x}{2}\right)}\right)}{\left(\frac{4x}{89} \cdot x + \frac{x}{2} \cdot x + \frac{4x}{89} \cdot \frac{x}{2}\right)^2 - 4 \cdot \frac{4x}{89} \cdot \frac{x}{2} \cdot x^2}$$

$$= \frac{\frac{2x^3}{89}\left(\frac{93x^2}{178} - \frac{18x^2}{89}\right)}{\left(\frac{101x^2}{178}\right)^2 - \frac{8x^4}{89}} = \frac{\frac{228x^5}{178^2}}{\frac{7353x^4}{178^2}} = \frac{4x}{129}$$

$$\therefore \quad r_4 = \frac{4x}{129} \approx 0.031007751937984496x \qquad \dots\dots\dots(11)$$

The above Eq(11) gives the radius of fourth inscribed circle (as shown in the above figure-7) in terms of edge length $x$ of square ABCD.

**Radius of fifth inscribed circle ($r_5$):** Similarly, the radius $r_5$ can be computed by substituting $a = r_4$, $b = x/2$ and $R = x$ in the above Eq(1) as follows

$$c_{min} = \frac{abR\left(aR + bR - ab - 2\sqrt{abR(R - a - b)}\right)}{(aR + bR + ab)^2 - 4abR^2}$$

$$\Rightarrow r_5 = \frac{\frac{4x}{129} \cdot \frac{x}{2} \cdot x \left(\frac{4x}{129} \cdot x + \frac{x}{2} \cdot x - \frac{4x}{129} \cdot \frac{x}{2} - 2\sqrt{\frac{4x}{129} \cdot \frac{x}{2} \cdot x \left(x - \frac{4x}{129} - \frac{x}{2}\right)}\right)}{\left(\frac{4x}{129} \cdot x + \frac{x}{2} \cdot x + \frac{4x}{129} \cdot \frac{x}{2}\right)^2 - 4 \cdot \frac{4x}{129} \cdot \frac{x}{2} \cdot x^2}$$

$$= \frac{\frac{2x^3}{129}\left(\frac{133x^2}{258} - \frac{22x^2}{129}\right)}{\left(\frac{141x^2}{258}\right)^2 - \frac{8x^4}{129}} = \frac{\frac{356x^5}{258^2}}{\frac{15753x^4}{258^2}} = \frac{4x}{177}$$

$$\therefore \quad r_5 = \frac{4x}{177} \approx 0.022598870056497175x \qquad \dots\dots\dots(12)$$

The above Eq(12) gives the radius of fifth inscribed circle (as shown in the above figure-7) in terms of edge length $x$ of square ABCD.

Similarly, the radii $r_6, r_7, r_8, \ldots r_n$ of all other inscribed circles can be found out.

It is worth noticing that the value of radius of each inscribed circle comes out in form of p/q whose numerator is always constant $4x$ ($x$ is constant side of square) but their denominators; 33, 57, 89, 129, 177, …….. form a series in which the differences between consecutive terms form an Arithmetic Progression (A.P.). Let this unknown series of the denominators be expressed in form of the sum of $n$ number of terms as follows

$$S_n = 33 + 57 + 89 + 129 + 177 + \dots\dots\dots\dots + T_n \qquad (1) \qquad (T_n \text{ is nth term})$$

$$S_n = \phantom{33 +} 33 + 57 + 89 + 129 + 177 + \dots\dots\dots\dots + T_n \qquad (2) \qquad (\text{Rewriting the series})$$

Now, subtracting the above series (2) from (1) (i.e. taking difference of corresponding terms of above two series), we obtain

$$0 = 33 + 24 + 32 + 40 + 48 + \ldots\ldots\ldots\ldots -T_n \qquad \text{(Total } (n+1) \text{ terms on RHS)}$$

$$T_n - 33 = 24 + 32 + 40 + 48 + \ldots\ldots\ldots\ldots \qquad \text{(Total } (n-1) \text{ terms on RHS)}$$

The series on the RHS is an Arithmetic Progression (A.P.) with first term $a = 24$, common difference $d = 8$ and total $(n-1)$ terms. Now, taking the sum of all (n-1) terms of above A.P. (on RHS) as follows

$$T_n - 33 = \frac{(n-1)}{2}(2 \cdot 24 + (n-1-1) \cdot 8) = 4(n-1)(n+4)$$

$$T_n = 33 + 4n^2 + 12n - 16$$

$$\boldsymbol{T_n = 4n^2 + 12n + 17}$$

Since the numerator of p/q form of radius of inscribed circle remains constant equal to $4x$ but the denominator of nth p/q form i.e. radius $r_n$ of nth inscribed circle is $T_n = 4n^2 + 12n + 17$. Therefore, following the symmetry/pattern in the denominators of p/q form of radii $r_1, r_2, r_3, r_4, r_5 \ldots \ldots r_n$ of inscribed circles (as derived in above section 4), the radius $r_n$ of nth inscribed circle is given as follows

$$r_n = \frac{4x}{T_n} = \frac{4x}{4n^2 + 12n + 17}$$

$$\therefore \ \boldsymbol{r_n = \frac{4x}{4n^2 + 12n + 17}} \qquad \ldots\ldots\ldots (12)$$

The above Eq(12) is the generalized formula to directly and analytically compute the radius $r_n$ of nth circle packed/inscribed in the plane region bounded by two semicircles and a quarter circle inside the square of given edge length $x$ (as shown in the above figure-7).

Now, substituting the different values of natural number $n = 1, 2, 3, 4 \ldots\ldots$ in the above generalized formula i.e. Eq(12), the radii of inscribed circles are obtained as follows

$$r_1 = \frac{4x}{33}, r_2 = \frac{4x}{57}, r_3 = \frac{4x}{89}, r_4 = \frac{4x}{129}, r_5 = \frac{4x}{177}, r_6 = \frac{4x}{233}, r_7 = \frac{4x}{297}, r_8 = \frac{4x}{369}, \ldots\ldots$$

The recurrence relation i.e. Eq(7) is equally applicable in this case to find the radius of any inscribed circle as follows

$$\boldsymbol{r_{n+1} = x\left(\frac{r_n^2 + xr_n - 2r_n\sqrt{xr_n - 2r_n^2}}{9r_n^2 - 2xr_n + x^2}\right)} \qquad \ldots\ldots\ldots (7)$$

Where, $\boldsymbol{r_1 = \frac{4x}{33}} \quad \forall \ x > 0, \ x \in R, \ n \in N \ \Rightarrow \ \boldsymbol{r_{n+1} = f(r_n)}$

However a simple MATLAB code has been developed using recurrence relation i.e. Eq(7) to compute the radius of all n number of circles packed/inscribed in the plane region bounded by two semicircles and a quarter circle inside a square of given edge length, as follows

## 4.1 MATLAB code-2

Using the above recurrence relation i.e. Eq(7), a for-loop can easily be used to develop a simple MATLAB code

```
clc;
clear all;
close all;
x=input('Enter the edge or side of square (in required unit) = ');
n=input('Enter the number of circles to be packed, n = ');

disp('Column matrix of radii of n no. of packed/inscribed circles, r(i)(units)=');
r(1)=(4*x)/33;
disp(r(1));
for i=1:1:n-1
    r(1)=(4*x)/33;
%HCR's recurrence relation to compute radii of inscribed circles
r(i+1)=double(x*((r(i))^2+x*r(i)-2*r(i)*sqrt(x*r(i)-2*(r(i))^2))/(9*(r(i))^2-
2*x*r(i)+x^2));
A=[r(i+1)];
disp(A);
end
```

Above MATLAB code of recurrence relation i.e. Eq(7) requires two input values 1) side $x$ of square and 2) number $n$ of circles to be inscribed/packed to compute the radii of all n number of circles packed or inscribed in the plane region bounded by two semi-circles and a quarter circle inside the square as shown in the figure-7 above.

## 4.2 Results

After entering the values of side of unit square i.e. $x = 1$ and number of circles to packed/inscribed i.e. $n = 100$, the following are the radii of all 100 packed/inscribed as tabulated below.

**Inputs**

```
Enter the edge or side of square (in required unit) =
1
Enter the number of circles to be packed, n =
100
```

**Outputs**

```
Column matrix of radii of n no. of packed/inscribed circles, r(i)(units)=
```

**Table-2: Ratio of radius of packed/inscribed circle to the side of square (n=100 circles)**

| | | | |
|---|---|---|---|
| 0.1212 | 0.0013 | 3.6255e-04 | 1.6644e-04 |
| 0.0702 | 0.0012 | 3.4913e-04 | 1.6223e-04 |
| 0.0449 | 0.0011 | 3.3645e-04 | 1.5817e-04 |
| 0.0310 | 0.0011 | 3.2444e-04 | 1.5427e-04 |
| 0.0226 | 0.0010 | 3.1306e-04 | 1.5051e-04 |
| 0.0172 | 9.4496e-04 | 3.0227e-04 | 1.4688e-04 |
| 0.0135 | 8.8948e-04 | 2.9203e-04 | 1.4338e-04 |
| 0.0108 | 8.3875e-04 | 2.8231e-04 | 1.4001e-04 |
| 0.0089 | 7.9224e-04 | 2.7306e-04 | 1.3676e-04 |
| 0.0074 | 7.4948e-04 | 2.6425e-04 | 1.3361e-04 |
| 0.0063 | 7.1010e-04 | 2.5587e-04 | 1.3058e-04 |
| 0.0054 | 6.7374e-04 | 2.4788e-04 | 1.2764e-04 |
| 0.0047 | 6.4010e-04 | 2.4025e-04 | 1.2481e-04 |
| 0.0041 | 6.0892e-04 | 2.3298e-04 | 1.2207e-04 |
| 0.0036 | 5.7996e-04 | 2.2603e-04 | 1.1941e-04 |
| 0.0032 | 5.5302e-04 | 2.1938e-04 | 1.1685e-04 |
| 0.0029 | 5.2791e-04 | 2.1303e-04 | 1.1436e-04 |
| 0.0026 | 5.0448e-04 | 2.0694e-04 | 1.1195e-04 |
| 0.0024 | 4.8257e-04 | 2.0112e-04 | 1.0962e-04 |
| 0.0022 | 4.6205e-04 | 1.9553e-04 | 1.0736e-04 |
| 0.0020 | 4.4282e-04 | 1.9018e-04 | 1.0517e-04 |
| 0.0018 | 4.2476e-04 | 1.8504e-04 | 1.0305e-04 |
| 0.0017 | 4.0779e-04 | 1.8011e-04 | 1.0099e-04 |
| 0.0015 | 3.9181e-04 | 1.7537e-04 | 9.8988e-05 |
| 0.0014 | 3.7675e-04 | 1.7082e-04 | 9.7047e-05 |

It's worth noticing that each value of above Table-2 shows the ratio of radius of a packed/inscribed circle to the side of square (i.e. taken as unity). Therefore, the exact value of radius of any of packed/inscribed circle of any order (i.e. up to $n = 100$) can be obtained by directly multiplying given side of square $x$ with the corresponding value of ratio as obtained in the above Table-2 (each column has 25 values of ratio) for 100 packed circles.

**Note:** The output values of radii of inscribed circles are obtained in form of a column from the above MALAB code-2 which have been rearranged in sequence in form of above Table-2 with four columns. All the values (i.e. obtained by recurrence relation) of Table-2 can be easily verified and validated by comparing them with the values obtained by above analytical formula i.e. Eq(12).

## 5. Packing of circles inside the sector of given radius and aperture angle

Consider a sector OAB consisting of two radii OA and OB each of length $R$ and an intercepted arc AB such that the central angle exerted by arc AB at the centre O is $\theta_1$ (as shown in the figure-10).

Now, consider $n$ number of circles with centres $C_1, C_2, C_3, \ldots\ldots C_n$ and radii $r_1, r_2, r_3, \ldots\ldots r_n$ respectively packed/inscribed in the sector OAB (as shown in the figure-10). Let's first find the radius $r_1$ of first packed/inscribed circle then the radii of all other inscribed circles $r_2, r_3, \ldots\ldots r_n$ can be easily computed by using recurrence relation.

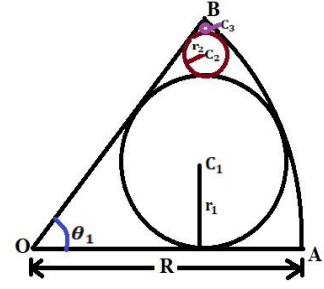

**Figure-10: Infinite no. of small circles with centers $C_1, C_2, C_3$ ... can be packed/inscribed in the sector OAB.**

Join the centres $C_1$ and $C_2$ and to the centre O of sector OAB. Drop the perpendiculars $C_1M$ and $C_2N$ from the centres $C_1$ and $C_2$ to the radii OA and OB respectively (as shown in the figure-11 below).

In right $\Delta C_1 MO$, We have, $OC_1 = R - r_1$, $C_1 M = r_1$ (see figure-11 below)

$$\sin \angle C_1 OM = \frac{C_1 M}{OC_1} \quad \Rightarrow \quad \sin \frac{\theta_1}{2} = \frac{r_1}{R - r_1}$$

$$R \sin \frac{\theta_1}{2} - r_1 \sin \frac{\theta_1}{2} = r_1 \quad \Rightarrow \quad r_1 \left(1 + \sin \frac{\theta_1}{2}\right) = R \sin \frac{\theta_1}{2}$$

$$\therefore \; \boldsymbol{r_1 = \frac{R \sin \frac{\theta_1}{2}}{1 + \sin \frac{\theta_1}{2}}} \qquad (\forall \; 0 < \theta_1 < \pi) \qquad \ldots\ldots\ldots (13)$$

The above Eq(13) gives the radius of first inscribed circle (as shown in the above figure-10) in terms of radius $R$ and central angle $\theta_1$ of sector OAB.

In right $\Delta C_2 NO$, We have, $C_2 N = r_2$, $OC_2 = R - r_2$.

Let $\angle BOD = \theta_2$ (see figure-11).

$$\Rightarrow \sin \angle C_2 ON = \frac{C_2 N}{OC_2} \quad \Rightarrow \quad \sin \frac{\theta_2}{2} = \frac{r_2}{R - r_2}$$

$$\sin \frac{\theta_2}{2} = \frac{r_2}{R\left(1 - \frac{r_2}{R}\right)} = \frac{\frac{r_2}{R}}{1 - \frac{r_2}{R}} = \boldsymbol{\frac{x}{1 - x}}$$

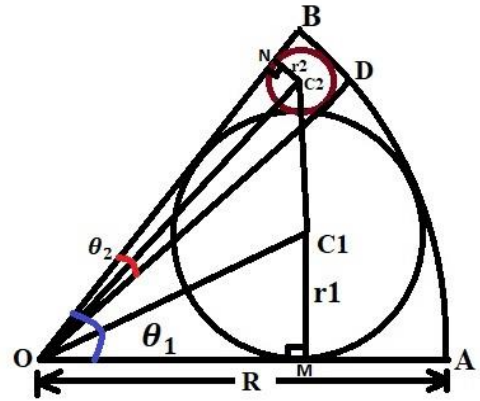

**Figure-11: Infinite no. of small circles with centers $C_1, C_2, C_3$ ... can be packed/inscribed in the sector OAB.**

Where, $x = \frac{r_2}{R}$ is some variable depending on the value of radius $r_2$ for given radius $R$ of sector AOB.

$$\Rightarrow \cos \frac{\theta_2}{2} = \sqrt{1 - \sin^2 \frac{\theta_2}{2}} \qquad \left(\because 0 < \theta_2 < \theta_1 < \pi \Leftrightarrow 0 < \frac{\theta_2}{2} < \frac{\theta_1}{2} < \frac{\pi}{2} \Rightarrow \cos \frac{\theta_2}{2} > 0\right)$$

$$\boldsymbol{\cos \frac{\theta_2}{2} = \sqrt{1 - \left(\frac{x}{1-x}\right)^2} = \frac{\sqrt{1 - 2x}}{1 - x}}$$

In $\Delta C_1 O C_2$, We have, $OC_1 = R - r_1$, $OC_2 = R - r_2$ and $C_1 C_2 = r_1 + r_2$. (see above figure-11).

$$\Rightarrow \angle C_1OC_2 = \angle C_1OB - \angle C_2OB = \frac{\angle AOB}{2} - \frac{\angle BOD}{2} = \frac{\theta_1}{2} - \frac{\theta_2}{2}$$

Now using **cosine rule** in $\Delta C_1OC_2$ as follows

$$\cos \angle C_1OC_2 = \frac{(OC_1)^2 + (OC_2)^2 - (C_1C_2)^2}{2(OC_1)(OC_2)}$$

$$\cos\left(\frac{\theta_1}{2} - \frac{\theta_2}{2}\right) = \frac{(R - r_1)^2 + (R - r_2)^2 - (r_1 + r_2)^2}{2(R - r_1)(R - r_2)}$$

$$\cos\frac{\theta_1}{2}\cos\frac{\theta_2}{2} + \sin\frac{\theta_1}{2}\sin\frac{\theta_2}{2} = \frac{R^2 - Rr_1 - Rr_2 - r_1r_2}{(R - r_1)(R - r_2)}$$

Snow, substituting the values of $\cos\frac{\theta_2}{2}$, $\sin\frac{\theta_2}{2}$ and $r_1$ in the above equation as follows

$$\cos\frac{\theta_1}{2} \cdot \frac{\sqrt{1 - 2x}}{1 - x} + \sin\frac{\theta_1}{2} \cdot \frac{x}{1 - x} = \frac{R^2 - R \cdot \frac{R\sin\frac{\theta_1}{2}}{1 + \sin\frac{\theta_1}{2}} - Rr_2 - \frac{R\sin\frac{\theta_1}{2}}{1 + \sin\frac{\theta_1}{2}} \cdot r_2}{\left(R - \frac{R\sin\frac{\theta_1}{2}}{1 + \sin\frac{\theta_1}{2}}\right)(R - r_2)}$$

$$\frac{\sqrt{1 - 2x}\cos\frac{\theta_1}{2} + x\sin\frac{\theta_1}{2}}{1 - x} = \frac{\frac{R^2}{1 + \sin\frac{\theta_1}{2}}\left(1 + \sin\frac{\theta_1}{2} - \sin\frac{\theta_1}{2} - \frac{r_2}{R} - \frac{r_2}{R}\sin\frac{\theta_1}{2} - \frac{r_2}{R}\sin\frac{\theta_1}{2}\right)}{\frac{R^2}{1 + \sin\frac{\theta_1}{2}}\left(1 + \sin\frac{\theta_1}{2} - \sin\frac{\theta_1}{2}\right)\left(1 - \frac{r_2}{R}\right)}$$

$$\frac{\sqrt{1 - 2x}\cos\frac{\theta_1}{2} + x\sin\frac{\theta_1}{2}}{1 - x} = \frac{\left(1 - x - 2x\sin\frac{\theta_1}{2}\right)}{1 - x} \qquad \left(\because \frac{r_2}{R} = x\right)$$

$$\frac{\sqrt{1 - 2x}\cos\frac{\theta_1}{2} + x\sin\frac{\theta_1}{2} - 1 + x + 2x\sin\frac{\theta_1}{2}}{1 - x} = 0$$

$$\sqrt{1 - 2x}\cos\frac{\theta_1}{2} + \left(1 + 3\sin\frac{\theta_1}{2}\right)x - 1 = 0 \qquad \left(\because x - 1 \neq 0 \Rightarrow x \neq 1 \vee r_2 \neq R\right)$$

$$\sqrt{1 - 2x}\cos\frac{\theta_1}{2} = 1 - \left(1 + 3\sin\frac{\theta_1}{2}\right)x$$

Taking the squares on both the (positive) sides of the above equation as follows

$$\left(\sqrt{1 - 2x}\cos\frac{\theta_1}{2}\right)^2 = \left(1 - \left(1 + 3\sin\frac{\theta_1}{2}\right)x\right)^2$$

$$(1 - 2x)\cos^2\frac{\theta_1}{2} = \left(1 + 3\sin\frac{\theta_1}{2}\right)^2 x^2 - 2\left(1 + 3\sin\frac{\theta_1}{2}\right)x + 1$$

$$\left(1 + 3\sin\frac{\theta_1}{2}\right)^2 x^2 - 2\left(1 + 3\sin\frac{\theta_1}{2} - \cos^2\frac{\theta_1}{2}\right)x + 1 - \cos^2\frac{\theta_1}{2} = 0$$

$$\left(1 + 3\sin\frac{\theta_1}{2}\right)^2 x^2 - 2\sin\frac{\theta_1}{2}\left(3 + \sin\frac{\theta_1}{2}\right)x + \sin^2\frac{\theta_1}{2} = 0$$

Solving the above quadratic equation for $x$ using quadratic formula as follows

$$x = \frac{-\left(-2\sin\frac{\theta_1}{2}\left(3+\sin\frac{\theta_1}{2}\right)\right) \pm \sqrt{\left(-2\sin\frac{\theta_1}{2}\left(3+\sin\frac{\theta_1}{2}\right)\right)^2 - 4\left(1+3\sin\frac{\theta_1}{2}\right)^2 \sin^2\frac{\theta_1}{2}}}{2\left(1+3\sin\frac{\theta_1}{2}\right)^2}$$

$$\frac{r_2}{R} = \frac{2\sin\frac{\theta_1}{2}\left(3+\sin\frac{\theta_1}{2}\right) \pm 2\sin\frac{\theta_1}{2}\sqrt{8-8\sin^2\frac{\theta_1}{2}}}{2\left(1+3\sin\frac{\theta_1}{2}\right)^2} \qquad \left(\because x = \frac{r_2}{R}\right)$$

$$r_2 = \frac{R\sin\frac{\theta_1}{2}\left(3+\sin\frac{\theta_1}{2}\right) \pm 2\sqrt{2}\sin\frac{\theta_1}{2}\cos\frac{\theta_1}{2}}{\left(1+3\sin\frac{\theta_1}{2}\right)^2}$$

$$r_2 = \frac{R\left(\sin^2\frac{\theta_1}{2} + 3\sin\frac{\theta_1}{2} \pm \sqrt{2}\sin\theta_1\right)}{\left(1+3\sin\frac{\theta_1}{2}\right)^2}$$

**Case-1**: Taking positive sign

$$r_2 = \frac{R\left(\sin^2\frac{\theta_1}{2} + 3\sin\frac{\theta_1}{2} + \sqrt{2}\sin\theta_1\right)}{\left(1+3\sin\frac{\theta_1}{2}\right)^2}$$

The values of radii obtained from the above expression of $r_2$ when used as a recurrence relation have been compared with the values of radii obtained from the drawings and measurements. It has been verified that the above expression of the radius $r_2$ gives erroneous values of radii of the packed/inscribed circles (i.e. $n$ number of circles) in the sector because the recurrence relation obtained from the above value of $r_2$ gives the increasing values of radii of inscribed circles with decreasing size which is not possible. Therefore, this case is discarded.

**Case-2**: Taking negative sign

$$r_2 = \frac{R\left(\sin^2\frac{\theta_1}{2} + 3\sin\frac{\theta_1}{2} - \sqrt{2}\sin\theta_1\right)}{\left(1+3\sin\frac{\theta_1}{2}\right)^2}$$

It has been verified that the above expression of $r_2$ when used as a recurrence relation gives the decreasing values of radii of inscribed circles (i.e. $n$ number of circles) which are exactly same as obtained from the drawings and measurements. Hence this case is valid and the above value of $r_2$ is acceptable.

$$r_2 = \frac{R\left(\sin^2\frac{\theta_1}{2} + 3\sin\frac{\theta_1}{2} - \sqrt{2}\sin\theta_1\right)}{\left(1+3\sin\frac{\theta_1}{2}\right)^2} \qquad \ldots\ldots\ldots\ldots\ldots (14)$$

The above Eq(14) gives the value of radius $r_2$ of second circle inscribed in the sector OAB (as shown in the figure-11 above).

Now, substituting the value of $r_2$ from Eq(14) in the expression of $\sin\frac{\theta_2}{2}$ (i.e. derived above) to obtain the value angle $\theta_2$ as follows

$$\sin\frac{\theta_2}{2} = \frac{r_2}{R - r_2} = \frac{\dfrac{R\left(\sin^2\frac{\theta_1}{2} + 3\sin\frac{\theta_1}{2} - \sqrt{2}\sin\theta_1\right)}{\left(1 + 3\sin\frac{\theta_1}{2}\right)^2}}{R - \dfrac{R\left(\sin^2\frac{\theta_1}{2} + 3\sin\frac{\theta_1}{2} - \sqrt{2}\sin\theta_1\right)}{\left(1 + 3\sin\frac{\theta_1}{2}\right)^2}} = \frac{\sin^2\frac{\theta_1}{2} + 3\sin\frac{\theta_1}{2} - \sqrt{2}\sin\theta_1}{1 + 8\sin^2\frac{\theta_1}{2} + 3\sin\frac{\theta_1}{2} + \sqrt{2}\sin\theta_1}$$

$$\theta_2 = 2\sin^{-1}\left(\frac{\sin^2\frac{\theta_1}{2} + 3\sin\frac{\theta_1}{2} - \sqrt{2}\sin\theta_1}{\sqrt{2}\sin\theta_1 + 8\sin^2\frac{\theta_1}{2} + 3\sin\frac{\theta_1}{2} + 1}\right) \quad\quad \ldots\ldots\ldots\ldots (15)$$

The above Eq(15) gives the value of angle $\theta_2$ exerted by the second inscribed circle at the centre O of the sector OAB (as shown in the figure-11 above).

Now, consider the new sector BOD of radius $R$ and central angle $\theta_2$ which inscribes the second circle of radius $r_2$ and centre $C_2$ and the third circle of radius $r_2$ and centre $C_2$ (as shown in the figure-12). This case becomes similar to the precious case in the above figure-11 in which the sector OAB of radius $R$ and central angle $\theta_1$ inscribes the circles with centres $C_1$ and $C_2$.

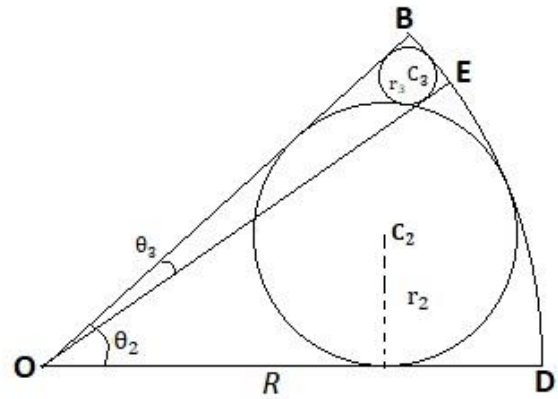

**Figure-12:** The sectors BOD and BOE are similar to the original sector OAB of radius $R$ and central angle $\theta_1$.

Therefore using the symmetry in the sectors OAB and OBD (see figure-11 & figure-12), the radius $r_3$ of third packed/inscribed circle with centre $C_3$ can be computed by substituting $\theta_1 = \theta_2$ in the above Eq(14) as follows

$$r_3 = \frac{R\left(\sin^2\frac{\theta_2}{2} + 3\sin\frac{\theta_2}{2} - \sqrt{2}\sin\theta_2\right)}{\left(1 + 3\sin\frac{\theta_2}{2}\right)^2}$$

The angle $\theta_3$ subtended by the third packed/inscribed circle with centre $C_3$ at the centre O of the sector OBD (i.e. sector OAB ) (see figure-12) can be computed by substituting $\theta_1 = \theta_2$ in the above Eq(15) as follows

$$\theta_3 = 2\sin^{-1}\left(\frac{\sin^2\frac{\theta_2}{2} + 3\sin\frac{\theta_2}{2} - \sqrt{2}\sin\theta_2}{\sqrt{2}\sin\theta_2 + 8\sin^2\frac{\theta_2}{2} + 3\sin\frac{\theta_2}{2} + 1}\right)$$

Similarly, consider the sector OBE of radius $R$ and central angle $\theta_3$ which inscribes the third circle of radius $r_3$ and centre $C_3$. The radius $r_4$ and the subtended angle $\theta_4$ of the fourth packed/inscribed circle are obtained by substituting $\theta_1 = \theta_3$ in the above Eq(14) and Eq(15) as follows

$$r_4 = \frac{R\left(\sin^2\frac{\theta_3}{2} + 3\sin\frac{\theta_3}{2} - \sqrt{2}\sin\theta_3\right)}{\left(1 + 3\sin\frac{\theta_3}{2}\right)^2}$$

$$\theta_4 = 2\sin^{-1}\left(\frac{\sin^2\frac{\theta_3}{2} + 3\sin\frac{\theta_3}{2} - \sqrt{2}\sin\theta_3}{\sqrt{2}\sin\theta_3 + 8\sin^2\frac{\theta_3}{2} + 3\sin\frac{\theta_3}{2} + 1}\right)$$

Similarly, the radii $r_5, r_6, r_7 \ldots$ and exerted angles $\theta_5, \theta_6, \theta_7 \ldots$ of all other circles packed/inscribed in the original sector OAB of radius $R$ and the cetral angle $\theta_1$ can be determined by substituting $\theta_1$ by the previous value in the above Eq(14) and Eq(15).

It can be concluded that the radius $r_{n+1}$ and subtended angle $\theta_{n+1}$ of say $(n+1)$th packed/inscribed circle depends on the previous value of angle $\theta_n$ i.e. angle exerted by $n$th packed/inscribed circle. Therefore, the radius $r_{n+1}$ and subtended angle $\theta_{n+1}$ of $(n+1)$th packed/inscribed circle are determined by substituting $\theta_1 = \theta_n$ in the above Eq(14) and Eq(15) as follows

$$r_{n+1} = \frac{R\left(\sin^2\frac{\theta_n}{2} + 3\sin\frac{\theta_n}{2} - \sqrt{2}\sin\theta_n\right)}{\left(1 + 3\sin\frac{\theta_n}{2}\right)^2} \quad \ldots \ldots \ldots (16)$$

$$\theta_{n+1} = 2\sin^{-1}\left(\frac{\sin^2\frac{\theta_n}{2} + 3\sin\frac{\theta_n}{2} - \sqrt{2}\sin\theta_n}{\sqrt{2}\sin\theta_n + 8\sin^2\frac{\theta_n}{2} + 3\sin\frac{\theta_n}{2} + 1}\right) \quad \ldots \ldots \ldots (17)$$

Where, $r_1 = \frac{R\sin\frac{\theta_1}{2}}{1+\sin\frac{\theta_1}{2}} \quad \forall\ 0 < \theta_1 < \pi,\ n \in N \quad \Rightarrow \quad r_{n+1} = f(\theta_n),\ \theta_{n+1} = g(\theta_n)$

The above recurrence relations i.e. Eq(16) and Eq(17) can be used to determine the radius and the angle subtended by any of infinite number of circles packed/inscribed in the sector of radius $R$ and central angle $\theta_1$ (refer to the above figure-10).

It is interesting to note that the radius and the subtended angle (i.e. obtained from Eq(16) & Eq(17)) of packed/inscribed circle decrease dramatically i.e. their decay is much faster than the exponential decay. In other words, the values of radius ($r_i$) and the subtended angle ($\theta_i$) approach to zero faster while packing a few number of circles in the sector of finite radius and central angle.

A simple MATLAB code has been developed using the above recurrence relations i.e. Eq(16) and Eq(17) to compute the radii and the subtended angles of all n number of circles packed/inscribed in the sector of given radius and central angle, as follows

### 5.1 MATLAB code-3

Using the above recurrence relations i.e. Eq(16) and Eq(17), a for-loop can easily be used to develop a simple MATLAB code

```
clc;
clear all;
close all;
R=input('Enter the radius of original sector, R (in required unit) = ');
t1=input('Enter the central angle of sector, theta_1 (in degree) = ');
n=input('Enter the number of circles to be packed, n = ');
r1=double((R*sind(t1/2))/(1+sind(t1/2)));
A=[r1 t1];
disp('(nx2) matrix of radii & angles subtended by inscribed circles, [r(units) theta(degrees)] =');
disp(A);
for i=1:n-1
    t(1)=(t1/2)*(pi/180);
```

```
t(i+1)=asin(((sin(t(i)))^2+3*sin(t(i))-
sqrt(2)*sin(2*t(i)))/(8*(sin(t(i)))^2+3*sin(t(i))+sqrt(2)*sin(2*t(i))+1));
r(i+1)=double((R*((sin(t(i)))^2+3*sin(t(i))-
sqrt(2)*sin(2*t(i))))/(1+3*sin(t(i)))^2);
B=[r(i+1) (360/pi)*t(i+1)];
disp(B);
end
```

### 5.2 Results

Using the above MATLAB code, the values of radii $r_i$ (in given unit) and the corresponding values of angles $\theta_i$ (in degree) subtended by 15 circles packed/inscribed in the sectors of unit radius ($R = 1$) and central angles $\theta_1 = 30º, 45º, 60º, 90º, 120º, 135º, 150º$, are obtained and arranged in the Table-3 below.

**Table-3: Radii and angles subtended by the circles packed/inscribed in a sector of unit radius.**

| $i$ | $r_i$ | $\theta_i$ | $r_i$ | $\theta_i$ | $r_i$ | $\theta_i$ | $r_i$ | $\theta_i$ | $r_i$ | $\theta_i$ | $r_i$ | $\theta_i$ | $r_i$ | $\theta_i$ |
|---|---|---|---|---|---|---|---|---|---|---|---|---|---|---|
| 1 | 0.2056 | 30.0000 | 0.2768 | 45.0000 | 0.3333 | 60.0000 | 0.4142 | 90.0000 | 0.4641 | 120.0000 | 0.4802 | 135.0000 | 0.4913 | 150.0000 |
| 2 | 0.0432 | 5.1759 | 0.0638 | 7.8185 | 0.0840 | 10.5288 | 0.1239 | 16.2602 | 0.1640 | 22.6285 | 0.1845 | 26.1597 | 0.2056 | 30.0000 |
| 3 | 0.0077 | 0.8882 | 0.0116 | 1.3420 | 0.0155 | 1.8077 | 0.0238 | 2.7944 | 0.0329 | 3.8947 | 0.0378 | 4.5073 | 0.0432 | 5.1759 |
| 4 | 0.0013 | 0.1524 | 0.0020 | 0.2302 | 0.0027 | 0.3102 | 0.0042 | 0.4795 | 0.0058 | 0.6683 | 0.0067 | 0.7734 | 0.0077 | 0.8882 |
| 5 | 0.0002 | 0.0261 | 0.0003 | 0.0395 | 0.0005 | 0.0532 | 0.0007 | 0.0823 | 0.0010 | 0.1147 | 0.0012 | 0.1327 | 0.0013 | 0.1524 |
| 6 | $39.08 \times 10^{-6}$ | 0.0045 | $59.14 \times 10^{-6}$ | 0.0068 | $79.65 \times 10^{-6}$ | 0.0091 | $12.32 \times 10^{-5}$ | 0.0141 | $17.17 \times 10^{-5}$ | 0.0197 | $19.86 \times 10^{-5}$ | 0.0228 | $22.81 \times 10^{-5}$ | 0.0261 |
| 7 | $67 \times 10^{-7}$ | $76.97 \times 10^{-5}$ | $10.18 \times 10^{-6}$ | 0.0012 | $13.62 \times 10^{-6}$ | 0.0016 | $21.11 \times 10^{-6}$ | 0.0024 | $29.49 \times 10^{-6}$ | 0.0034 | $34.14 \times 10^{-6}$ | 0.0039 | $39.08 \times 10^{-6}$ | 0.0045 |
| 8 | $12 \times 10^{-7}$ | $13.21 \times 10^{-5}$ | $17 \times 10^{-7}$ | $19.95 \times 10^{-5}$ | $23 \times 10^{-7}$ | $26.88 \times 10^{-5}$ | $36 \times 10^{-7}$ | $41.55 \times 10^{-5}$ | $51 \times 10^{-7}$ | $57.91 \times 10^{-5}$ | $58 \times 10^{-7}$ | $67.02 \times 10^{-5}$ | $67 \times 10^{-7}$ | $76.97 \times 10^{-5}$ |
| 9 | $20 \times 10^{-8}$ | $22.66 \times 10^{-6}$ | $30 \times 10^{-8}$ | $34.23 \times 10^{-6}$ | $40 \times 10^{-8}$ | $46.11 \times 10^{-6}$ | $62 \times 10^{-8}$ | $71.28 \times 10^{-6}$ | $87 \times 10^{-8}$ | $99.36 \times 10^{-6}$ | $10 \times 10^{-7}$ | $11.5 \times 10^{-5}$ | $12 \times 10^{-7}$ | $13.21 \times 10^{-5}$ |
| 10 | $34 \times 10^{-9}$ | $38.87 \times 10^{-7}$ | $51 \times 10^{-9}$ | $58.73 \times 10^{-7}$ | $69 \times 10^{-9}$ | $79.12 \times 10^{-7}$ | $11 \times 10^{-8}$ | $12.23 \times 10^{-6}$ | $15 \times 10^{-8}$ | $17.05 \times 10^{-6}$ | $17 \times 10^{-8}$ | $19.73 \times 10^{-6}$ | $20 \times 10^{-8}$ | $22.66 \times 10^{-6}$ |
| 11 | $58 \times 10^{-10}$ | $66.7 \times 10^{-8}$ | $9 \times 10^{-9}$ | $10.08 \times 10^{-7}$ | $12 \times 10^{-9}$ | $13.57 \times 10^{-7}$ | $18 \times 10^{-9}$ | $20.98 \times 10^{-7}$ | $26 \times 10^{-9}$ | $29.25 \times 10^{-7}$ | $30 \times 10^{-9}$ | $33.85 \times 10^{-7}$ | $34 \times 10^{-9}$ | $38.87 \times 10^{-7}$ |
| 12 | $10 \times 10^{-10}$ | $11.44 \times 10^{-8}$ | $15 \times 10^{-10}$ | $17.29 \times 10^{-8}$ | $20 \times 10^{-10}$ | $23.29 \times 10^{-8}$ | $31 \times 10^{-10}$ | $36 \times 10^{-8}$ | $44 \times 10^{-10}$ | $50.18 \times 10^{-8}$ | $51 \times 10^{-10}$ | $58.08 \times 10^{-8}$ | $58 \times 10^{-10}$ | $66.7 \times 10^{-8}$ |
| 13 | $17 \times 10^{-11}$ | $19.63 \times 10^{-9}$ | $26 \times 10^{-11}$ | $29.66 \times 10^{-9}$ | $35 \times 10^{-11}$ | $39.96 \times 10^{-9}$ | $54 \times 10^{-11}$ | $61.77 \times 10^{-9}$ | $75 \times 10^{-11}$ | $86.1 \times 10^{-9}$ | $87 \times 10^{-11}$ | $99.65 \times 10^{-9}$ | $10 \times 10^{-10}$ | $11.44 \times 10^{-8}$ |
| 14 | $29 \times 10^{-12}$ | $33.69 \times 10^{-10}$ | $44 \times 10^{-12}$ | $50.9 \times 10^{-10}$ | $60 \times 10^{-12}$ | $68.56 \times 10^{-10}$ | $9 \times 10^{-11}$ | $10.6 \times 10^{-9}$ | $13 \times 10^{-11}$ | $14.77 \times 10^{-9}$ | $15 \times 10^{-11}$ | $17.1 \times 10^{-9}$ | $17 \times 10^{-11}$ | $19.63 \times 10^{-9}$ |
| 15 | $50 \times 10^{-13}$ | $57.80 \times 10^{-11}$ | $76 \times 10^{-13}$ | $87.32 \times 10^{-11}$ | $10 \times 10^{-12}$ | $11.76 \times 10^{-10}$ | $16 \times 10^{-12}$ | $18.18 \times 10^{-10}$ | $22 \times 10^{-12}$ | $25.35 \times 10^{-10}$ | $26 \times 10^{-12}$ | $29.33 \times 10^{-10}$ | $29 \times 10^{-12}$ | $33.69 \times 10^{-10}$ |

# 6. Packing of circles in the plane bounded by two identical circles & their common tangent

Consider two externally touching identical circles with centres A & B, and each of radius $R$ and their common tangent PQ. (as shown in the figure-13).

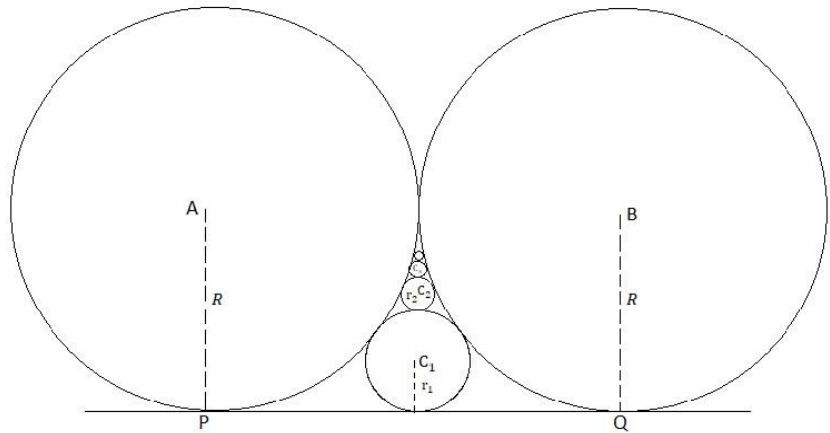

Figure-13: $n$ no. of circles are inscribed in the plane region bounded by two externally touching identical circles with centers A & B and each of radius $R$, and their common tangent PQ. The common tangent PQ can be assumed to be a circular arc of infinite radius.

Now, consider $n$ number of circles with centres $C_1, C_2, C_3, \ldots\ldots C_n$ and radii $r_1, r_2, r_3, \ldots\ldots r_n$ respectively packed/inscribed in the plane region bounded by two externally touching identical circles and their common tangent (as shown in the figure-13 above). Let's first find the radius $r_1$ of first packed/inscribed circle then the radii of all other inscribed circles $r_2, r_3, \ldots\ldots r_n$ can be easily computed by using generalized formula.

We know that the radius $(r)$ of inscribed circle which externally touches three given circles with radii say $a, b$ & $c$ touching each other externally (as shown in figure-14) is given by the following generalized formula [1]

$$r = \frac{abc}{2\sqrt{abc(a+b+c)} + (ab+bc+ca)} \qquad \ldots\ldots\ldots (18)$$

It's worth noticing that the common tangent PQ can be assumed to be the arc of a circle of infinite radius which externally touches two externally touching equal circles with centres A and B (see the above figure-13). Thus this case (figure-13) becomes similar to the case (figure-14) of three externally touching circles with radii $a = R, b = R$ and $c \to \infty$.

Now, the first circle with radius $r_1$ and centre $C_1$ is packed/inscribed in the plane region bounded two equal circles and an infinite circle (i.e. common tangent PQ). The first packed/inscribed circle with centre $C_1$ and radius $r_1$ externally touches three externally touching circles with radii $a = R, b = R$ and $c \to \infty$ (see the above figure-13).

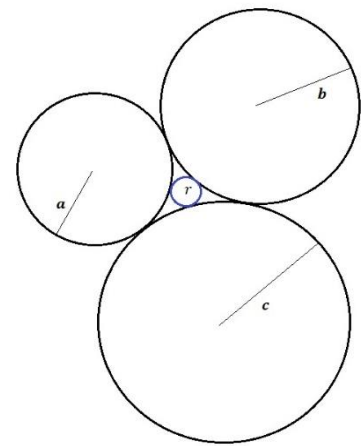

Figure-14: Three externally touching circles with radii $a, b,$ and $c$ inscribe a blue circle with radius $r$.

Therefore, the radius $r_1$ of first packed/inscribed circle (see figure-13 above) is determined by substituting $a = R$, $b = R$ in the above Eq(18) and taking its limit as $c \to \infty$ as follows

$$r_1 = \lim_{c \to \infty} \frac{R \cdot R \cdot c}{2\sqrt{R \cdot R \cdot c(R + R + c)} + (R \cdot R + R \cdot c + c \cdot R)} = \lim_{c \to \infty} \frac{R^2 c}{Rc\left(2\sqrt{\left(\frac{R}{c} + \frac{R}{c} + 1\right)} + \left(\frac{R}{c} + 2\right)\right)}$$

$$= \lim_{c \to \infty} \frac{R}{2\sqrt{\left(\frac{R}{c} + \frac{R}{c} + 1\right)} + \frac{R}{c} + 2} = \frac{R}{2\sqrt{(0 + 0 + 1)} + 0 + 2} = \frac{R}{4}$$

$$\Rightarrow \; r_1 = \frac{R}{4}$$

Now, consider the second packed/inscribed circle with radius $r_2$ and centre $C_2$ which externally touches the three externally touching circles with centres $A, B$ and $C_1$ and radii $R, R$ and $r_1$ respectively (see the above figure-13). Therefore, the radius $r_2$ is determined by substituting $a = R$, $b = R$ and $c = r_1 = \frac{R}{4}$ in the above Eq(18) as follows

$$r_2 = \frac{R \cdot R \cdot \frac{R}{4}}{2\sqrt{R \cdot R \cdot \frac{R}{4}\left(R + R + \frac{R}{4}\right)} + \left(R \cdot R + R \cdot \frac{R}{4} + \frac{R}{4} \cdot R\right)} = \frac{\frac{R^3}{4}}{\frac{3R^2}{2} + \frac{3R^2}{2}} = \frac{R}{12}$$

$$\Rightarrow \; r_2 = \frac{R}{12}$$

Similarly, consider the third packed/inscribed circle with radius $r_3$ and centre $C_3$ which externally touches the three externally touching circles with centres $A, B$ and $C_2$ and radii $R, R$ and $r_2$ respectively (see the above figure-13). Therefore, the radius $r_3$ is determined by substituting $a = R$, $b = R$ and $c = r_2 = \frac{R}{12}$ in the above Eq(18) as follows

$$r_3 = \frac{R \cdot R \cdot \frac{R}{12}}{2\sqrt{R \cdot R \cdot \frac{R}{12}\left(R + R + \frac{R}{12}\right)} + \left(R \cdot R + R \cdot \frac{R}{12} + \frac{R}{12} \cdot R\right)} = \frac{\frac{R^3}{12}}{\frac{5R^2}{6} + \frac{7R^2}{6}} = \frac{R}{24}$$

$$\Rightarrow \; r_3 = \frac{R}{24}$$

Similarly, the radius $r_4$ is determined by substituting $a = R$, $b = R$ and $c = r_3 = \frac{R}{24}$ in the above Eq(18) as follows

$$r_4 = \frac{R \cdot R \cdot \frac{R}{24}}{2\sqrt{R \cdot R \cdot \frac{R}{24}\left(R + R + \frac{R}{24}\right)} + \left(R \cdot R + R \cdot \frac{R}{24} + \frac{R}{24} \cdot R\right)} = \frac{\frac{R^3}{24}}{\frac{7R^2}{12} + \frac{13R^2}{12}} = \frac{R}{40}$$

$$\Rightarrow \; r_4 = \frac{R}{40}$$

Similarly, the radius $r_5$ is determined by substituting $a = R$, $b = R$ and $c = r_4 = \frac{R}{40}$ in the above Eq(18) as follows

$$r_5 = \frac{R \cdot R \cdot \frac{R}{40}}{2\sqrt{R \cdot R \cdot \frac{R}{40}\left(R + R + \frac{R}{40}\right)} + \left(R \cdot R + R \cdot \frac{R}{40} + \frac{R}{40} \cdot R\right)} = \frac{\frac{R^3}{40}}{\frac{9R^2}{20} + \frac{21R^2}{20}} = \frac{R}{60}$$

$$\Rightarrow \; r_5 = \frac{R}{60}$$

Similarly, the radius $r_6$ is determined by substituting $a = R$, $b = R$ and $c = r_5 = \frac{R}{60}$ in the above Eq(18) as follows

$$r_6 = \frac{R \cdot R \cdot \frac{R}{60}}{2\sqrt{R \cdot R \cdot \frac{R}{60}\left(R + R + \frac{R}{60}\right) + \left(R \cdot R + R \cdot \frac{R}{60} + \frac{R}{60} \cdot R\right)}} = \frac{\frac{R^3}{60}}{\frac{11R^2}{30} + \frac{31R^2}{30}} = \frac{R}{84}$$

$$\Rightarrow r_6 = \frac{R}{84}$$

Similarly, all other radii $r_7, r_8, r_9, \ldots r_n$ can be determined using the generalized formula i.e. the above Eq(18).

It is worth noticing that the value of radius of each inscribed circle comes out in form of $\frac{p}{q}$ whose numerator is always constant $R$ (i.e. given radius of equal circles) but their denominators; 4, 12, 24, 40, 60, 84, …….. form a series in which the differences between consecutive terms form an Arithmetic Progression (A.P.). Let this unknown series of the denominators be expressed in form of the sum of $n$ number of terms as follows

$$S_n = 4 + 12 + 24 + 40 + 60 + 84 \ldots\ldots\ldots\ldots. + T_n \qquad (1) \qquad (T_n \text{ is nth term})$$

$$S_n = \quad\ 4 + 12 + 24 + 40 + 60 + 84 \ldots\ldots\ldots\ldots. + T_n \qquad (2) \qquad (\text{Rewriting the series})$$

Now, subtracting the above series (2) from (1) (i.e. taking difference of corresponding terms of above two series), we obtain

$$0 = 4 + 8 + 12 + 16 + 20 + 24 + \ldots\ldots\ldots\ldots\ldots. - T_n \qquad (\text{Total } (n+1) \text{ terms on RHS})$$

$$T_n = 4 + 8 + 12 + 16 + 20 + 24 + \ldots\ldots\ldots\ldots\ldots. \qquad (\text{Total n terms on RHS})$$

The series on the RHS is an Arithmetic Progression (A.P.) with first term $a = 4$, common difference $d = 4$ and total $n$ terms. Now, taking the sum of all $n$ no. of terms of above A.P. (on RHS) as follows

$$T_n = \frac{n}{2}(2 \cdot 4 + (n-1) \cdot 4)$$

$$\boldsymbol{T_n = 2n(n+1)}$$

Since the numerator of $\frac{p}{q}$ form of radius of inscribed circle remains constant equal to $R$ but the denominator of nth $\frac{p}{q}$ form i.e. radius $r_n$ of nth inscribed circle is $T_n = 2n(n+1)$. Therefore, following the symmetry/pattern in the denominators of $\frac{p}{q}$ form of radii $r_1, r_2, r_3, r_4, r_5 \ldots\ldots r_n$, the radius $r_n$ of nth packed/inscribed circle is given as follows

$$r_n = \frac{R}{T_n} = \frac{R}{2n(n+1)}$$

$$\therefore\ \boldsymbol{r_n = \frac{R}{2n(n+1)}} \qquad \ldots\ldots\ldots (19)$$

The above Eq(19) is the generalized formula to directly and analytically compute the radius $r_n$ of nth circle packed/inscribed in the plane region bounded by two externally touching identical circles each of radius $R$ and their common tangent (as shown in the above figure-13).

Now, substituting the different values of natural number $n = 1, 2, 3, 4 \ldots\ldots$, in the above generalized formula i.e. Eq(19), the radii of inscribed circles are obtained as follows

$$r_1 = \frac{R}{4}, r_2 = \frac{R}{12}, r_3 = \frac{R}{24}, r_4 = \frac{R}{40}, r_5 = \frac{R}{60}, r_6 = \frac{R}{84}, r_7 = \frac{R}{112}, r_8 = \frac{R}{144}, \ldots\ldots$$

# 7. Closest 2D packing of identical circles in a regular hexagon

The closest, densest or hexagonal packing of circles over a plane is the closest/densest possible 2D arrangement of circles without any overlap such that each circle touches exactly six circles in its closest neighbourhood. (as shown in the figure-15). We are to derive the generalized formula to compute total number of identical circles most densely packed in a regular hexagon, the radius of smallest circle and the minimum length of string which encloses the identical circles of a given (finite) radius which are the most densely packed in a regular hexagon of given side.

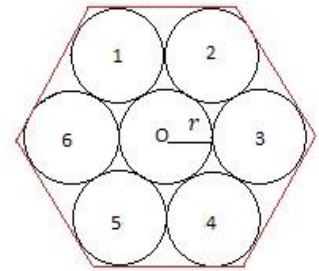

**Figure-15: The circle with center O and radius $r$ is tangent to 6 identical circles the most densely packed in a regular hexagon.**

Let's first define the following known parameters

$n$ = number of identical circles tangent to any of six sides of a regular hexagon in which the circles are the most densely packed

$r$ = radius of each of identical circles the most densely packed in the regular hexagon

The following important unknown parameters are to be determined in terms of known parameters ($n, r$)

$N$ = total number of identical circles which are the most densely packed in the regular hexagon

$N_v$ = total number of voids (i.e. bounded regions unoccupied by circles) in the regular hexagon

$a$ = side of regular hexagon

$R$ = radius of the smallest circle circumscribing all the circles the most densely packed in the regular hexagon

$L$ = minimum length of one turn of string (i.e. thin, perfectly flexible, and inextensible) wrapped around all the packed circles

Now, consider the following cases of the value of $n$ i.e. number of identical circles tangent to each side of regular hexagon as follows

**Case-1: $n = 2$ (i.e. two identical circles are tangent to each side of regular hexagon)**

Consider total $N$ number of identical circles each of radius $r$ the most densely packed in a regular hexagon ABCDEF of side $a$ such that two identical circles are tangent to each side say side ED (as shown in the figure-16(a) below).

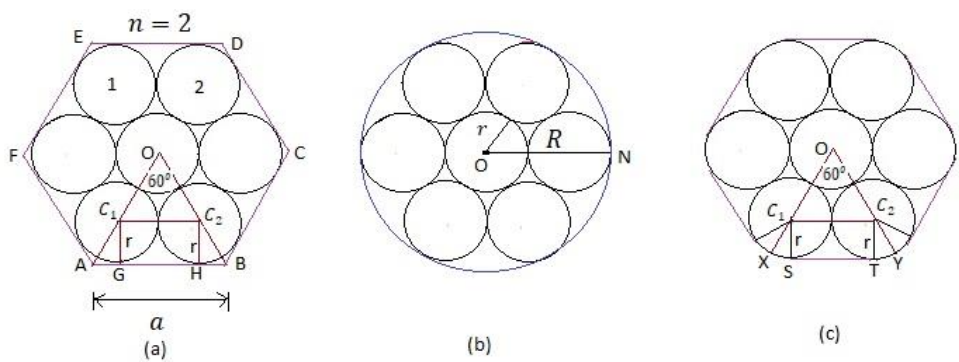

**Figure-16: (a) Two identical circles labeled with 1, 2 are tangent to (each) side ED of regular hexagon ABCDEF (n=2) (b) smallest circle circumscribing all the packed circles (c) minimum length of one turn of flexible & inextensible string wrapped around the identical circles most densely packed in a regular hexagon.**

In this case, unknown parameters $N, N_v, a, R,$ and $L$ (as defined above) can be easily determined as given below.

In order to find out the number of identical circles and number of voids, the centres of circles lying on the sides & vertices of internal regular hexagons are counted in CW circular fashion (as shown in the figure-17) as follows

$N$ = Total number of identical circles packed in regular hexagon
$$= 1 + 6$$

It is worth noticing that 6 voids are always present at 6 vertices/corners of the outermost regular hexagon ABCDEF (see figure-17).

$N_v$ = Total number of voids (unoccupied regions) in regular hexagon
$$= 6 + 2(6)$$

**Figure-17:** The centers of circles lie on inner hexagon which can be easily counted in CW circular fashion around the central circle

In right $\Delta C_1 GA$, We have $\angle AC_1G = \frac{\angle AOB}{2} = \frac{60^o}{2} = 30^o$, $C_1G = r$ (see the figure-16(a) above)

$$\tan \angle AC_1G = \frac{AG}{C_1G} \Rightarrow AG = C_1G \tan \angle AC_1G = r \tan 30^o = \frac{r}{\sqrt{3}}$$

$a$ = Side of regular hexagon ABCDEF = $AB = AG + GH + HB = GH + 2(AG) = 2r + \frac{2r}{\sqrt{3}} = (2-1)(2r) + \frac{2r}{\sqrt{3}}$

$R$ = Radius of smallest circumscribing circle = $ON = r + 2r = r + (2-1)(2r)$ (see figure 16(b) above)

In the figure-16(c) above, We have $\angle XC_1S = \frac{\angle C_1OC_2}{2} = \frac{60^o}{2} = \frac{\pi}{6}$, $C_1S = C_1X = r$

$$\Rightarrow \text{Arc } SX = C_1X\,(\angle XC_1S) = r\left(\frac{\pi}{6}\right) = \frac{\pi r}{6}$$

$\therefore XY = \text{arc } XS + ST + \text{arc } TY = ST + 2(\text{arc } XS) = 2r + 2\left(\frac{\pi r}{6}\right) = 2r + \frac{\pi r}{3}$ (see figure 16(c) above)

$L$ = Minimum length of one turn of string $= 6(XY) = 6\left(2r + \frac{\pi r}{3}\right) = 6\left((2-1)(2r) + \frac{\pi r}{3}\right)$

**Case-2: $n = 3$ (i.e. three identical circles are tangent to each side of regular hexagon)**

Similarly, consider total $N$ number of identical circles each of radius $r$ the most densely packed in a regular hexagon ABCDEF of side $a$ such that three identical circles are tangent to each side say side ED (as shown in the figure-18(a) below).

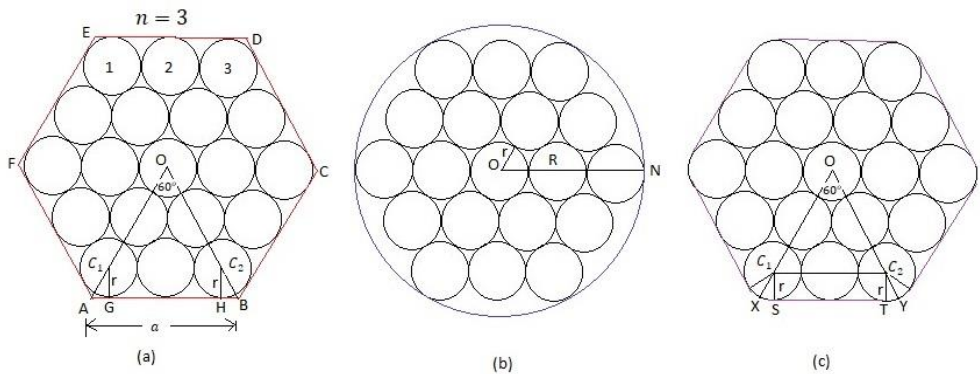

**Figure-18:** (a) Three identical circles labeled with 1, 2, 3 are tangent to (each) side ED of regular hexagon ABCDEF (n=3) (b) smallest circle circumscribing all the packed circles (c) minimum length of one turn of flexible & inextensible string wrapped around the identical circles most densely packed in a regular hexagon.

In this case, the unknown parameters $N, N_v, a, R,$ and $L$ can be easily determined as given below.

In order to find out the number of identical circles and number of voids, the centres of circles lying on the sides & vertices of internal regular hexagons are counted in CW circular fashion (as shown in the figure-19) as follows

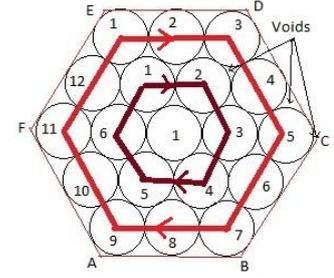

$N$ = Total number of identical circles packed in regular hexagon
$$= 1 + 6 + 12$$

It is worth noticing that 6 voids are always present at 6 vertices/corners of the outermost regular hexagon ABCDEF (see figure-19).

$N_v$ = Total number of voids (unoccupied regions) in regular hexagon
$$= 6 + 2(6 + 12)$$

**Figure-19: The centers of circles lie on inner hexagons which can be easily counted in CW circular fashion around the central circle.**

In right $\Delta C_1GA$, We have $\angle AC_1G = \frac{\angle AOB}{2} = \frac{60°}{2} = 30°$, $C_1G = r$ (see the figure-18(a) above)

$$\tan \angle AC_1G = \frac{AG}{C_1G} \Rightarrow AG = C_1G \tan \angle AC_1G = r \tan 30° = \frac{r}{\sqrt{3}}$$

$a$ = Side of regular hexagon ABCDEF = $AB = AG + GH + HB = GH + 2(AG) = 2(2r) + \frac{2r}{\sqrt{3}} = (3-1)(2r) + \frac{2r}{\sqrt{3}}$

$R$ = Radius of smallest circumscribing circle = $ON = r + 2(2r) = r + (3-1)(2r)$    (see figure 18(b) above)

In the figure-18(c) above, We have $\angle XC_1S = \frac{\angle C_1OC_2}{2} = \frac{60°}{2} = \frac{\pi}{6}$, $C_1S = C_1X = r$

$$\Rightarrow \text{Arc } SX = C_1X (\angle XC_1S) = r\left(\frac{\pi}{6}\right) = \frac{\pi r}{6}$$

$\therefore XY = \text{arc } XS + ST + \text{arc } TY = ST + 2(\text{arc } XS) = 4r + 2\left(\frac{\pi r}{6}\right) = 2(2r) + \frac{\pi r}{3}$    (see figure 18(c) above)

$L$ = Minimum length of one turn of string $= 6(XY) = 6\left(2(2r) + \frac{\pi r}{3}\right) = 6\left((3-1)(2r) + \frac{\pi r}{3}\right)$

**Case-3: $n = 4$ (i.e. four identical circles are tangent to each side of regular hexagon)**

Similarly, consider total $N$ number of identical circles each of radius $r$ the most densely packed in a regular hexagon ABCDEF of side $a$ such that four identical circles are tangent to each side say side ED (as shown in the figure-20(a) below).

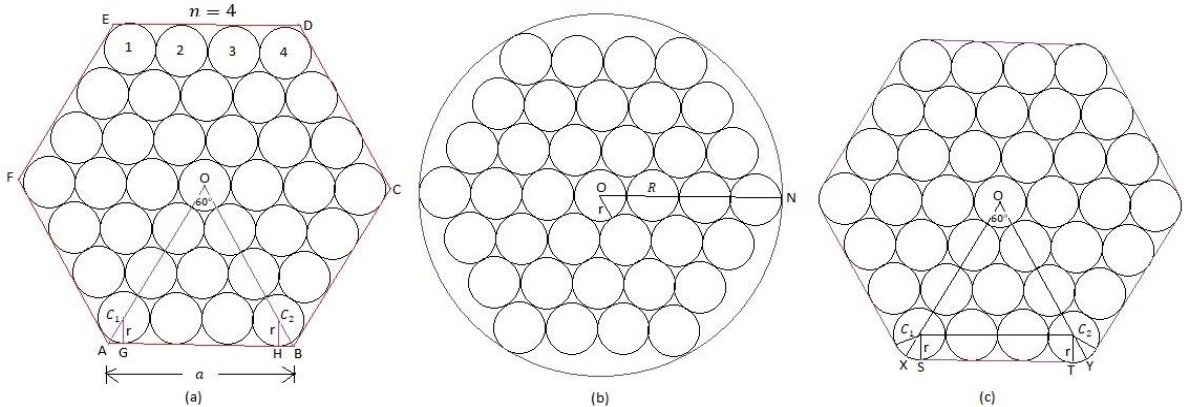

**Figure-20: (a) Four identical circles labeled with 1, 2, 3, 4 are tangent to (each) side ED of regular hexagon ABCDEF (n=4) (b) smallest circle circumscribing all the packed circles (c) minimum length of one turn of flexible & inextensible string wrapped around the identical circles most densely packed in a regular hexagon.**

In this case, the unknown parameters $N, N_v, a, R,$ and $L$ can be easily determined as given below.

In order to find out the number of identical circles and number of voids, the centres of circles lying on the sides & vertices of internal regular hexagons are counted in CW circular fashion (as shown in the figure-21) as follows

$N$ = Total number of identical circles packed in regular hexagon
$= 1 + 6 + 12 + 18$

It is worth noticing that 6 voids are always present at 6 vertices/corners of the outermost regular hexagon ABCDEF (see the figure-21).

$N_v$ = Total number of voids (unoccupied regions) in regular hexagon
$= 6 + 2(6 + 12 + 18)$

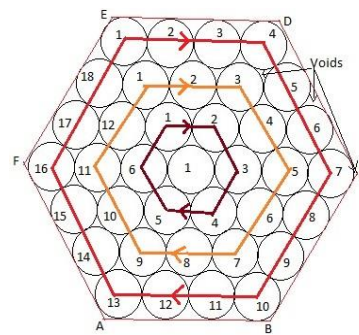

**Figure-21:** The centers of circles lie on inner hexagons which can be easily counted in CW circular fashion around the central circle.

In right $\Delta C_1 GA$, We have $\angle AC_1G = \frac{\angle AOB}{2} = \frac{60^o}{2} = 30^o$, $C_1G = r$ (see the figure-20(a) above)

$$\tan \angle AC_1G = \frac{AG}{C_1G} \Rightarrow AG = C_1G \tan \angle AC_1G = r \tan 30^o = \frac{r}{\sqrt{3}}$$

$a$ = Side of regular hexagon ABCDEF = $AB = AG + GH + HB = GH + 2(AG) = 3(2r) + \frac{2r}{\sqrt{3}} = (4-1)(2r) + \frac{2r}{\sqrt{3}}$

$R$ = Radius of smallest circumscribing circle = $ON = r + 3(2r) = r + (4-1)(2r)$ (see figure 20(b) above)

In the figure-20(c) above, We have $\angle XC_1S = \frac{\angle C_1OC_2}{2} = \frac{60^o}{2} = \frac{\pi}{6}$, $C_1S = C_1X = r$

$$\Rightarrow \text{Arc } SX = C_1X \ (\angle XC_1S) = r\left(\frac{\pi}{6}\right) = \frac{\pi r}{6}$$

$\therefore XY = \text{arc } XS + ST + \text{arc } TY = ST + 2(\text{arc } XS) = 6r + 2\left(\frac{\pi r}{6}\right) = 3(2r) + \frac{\pi r}{3}$ (see figure 20(c) above)

$L$ = Minimum length of one turn of string = $6(XY) = 6\left(3(2r) + \frac{\pi r}{3}\right) = 6\left((4-1)(2r) + \frac{\pi r}{3}\right)$

**Case-4: $n = 5$ (i.e. five identical circles are tangent to each side of regular hexagon)**

Similarly, consider total $N$ number of identical circles each of radius $r$ the most densely packed in a regular hexagon ABCDEF of side $a$ such that five identical circles are tangent to each side say side ED (as shown in the figure-23(a) below).

In this case, the unknown parameters $N, N_v, a, R,$ and $L$ can be easily determined as given below.

In order to find out the number of identical circles and number of voids, the centres of circles lying on the sides & vertices of internal regular hexagons are counted in CW circular fashion (as shown in the figure-22) as follows

$N$ = Total number of identical circles packed in regular hexagon
$= 1 + 6 + 12 + 18 + 24$

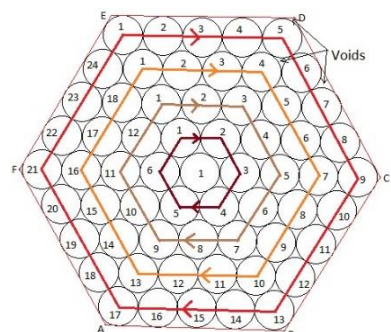

**Figure-22:** The centers of circles lie on inner hexagons which can be easily counted in CW circular fashion around the central circle.

It is worth noticing that 6 voids are always present at 6 vertices/corners of the outermost regular hexagon ABCDEF (see figure-22 above).

$$N_v = \text{Total number of voids (unoccupied regions) in regular hexagon}$$
$$= 6 + 2(6 + 12 + 18 + 24)$$

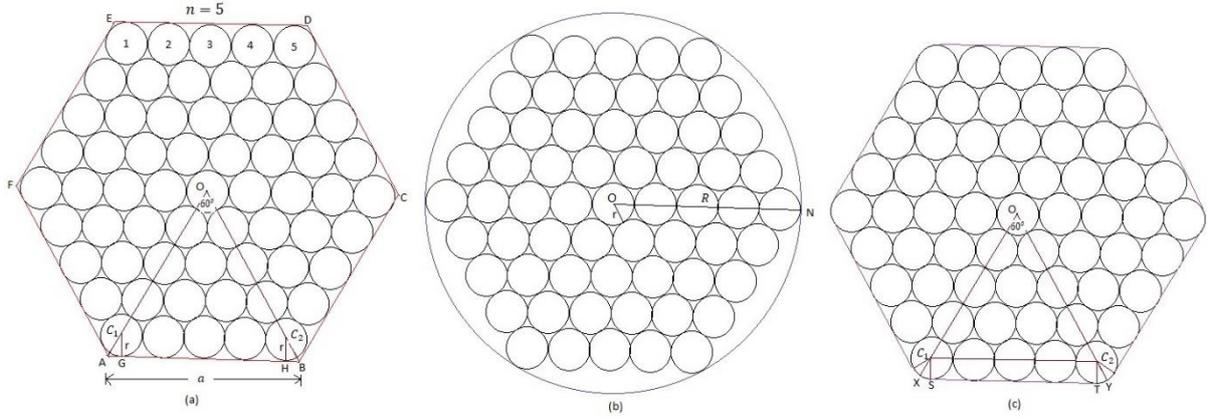

**Figure-23: (a) Five identical circles labeled with 1, 2, 3, 4, 5 are tangent to (each) side ED of regular hexagon ABCDEF (n=5) (b) smallest circle circumscribing all the packed circles (c) minimum length of one turn of flexible & inextensible string wrapped around the identical circles most densely packed in a regular hexagon.**

In right $\Delta C_1 GA$, We have $\angle AC_1 G = \frac{\angle AOB}{2} = \frac{60^o}{2} = 30^o$, $C_1 G = r$ (see the figure-23(a) above)

$\tan \angle AC_1 G = \frac{AG}{C_1 G} \Rightarrow AG = C_1 G \tan \angle AC_1 G = r \tan 30^o = \frac{r}{\sqrt{3}}$

$a = \text{Side of regular hexagon ABCDEF} = AB = AG + GH + HB = GH + 2(AG) = 4(2r) + \frac{2r}{\sqrt{3}} = (5-1)(2r) + \frac{2r}{\sqrt{3}}$

$R = \text{Radius of smallest circumscribing circle} = ON = r + 4(2r) = r + (5-1)(2r)$      (see figure 23(b) above)

In the figure-23(c) above, We have $\angle XC_1 S = \frac{\angle C_1 OC_2}{2} = \frac{60^o}{2} = \frac{\pi}{6}$, $C_1 S = C_1 X = r$

$$\Rightarrow \text{Arc } SX = C_1 X (\angle XC_1 S) = r \left(\frac{\pi}{6}\right) = \frac{\pi r}{6}$$

$\therefore XY = \text{arc } XS + ST + \text{arc } TY = ST + 2(\text{arc } XS) = 8r + 2\left(\frac{\pi r}{6}\right) = 4(2r) + \frac{\pi r}{3}$     (see figure 23(c) above)

$L = \text{Minimum length of one turn of string} = 6(XY) = 6\left(4(2r) + \frac{\pi r}{3}\right) = 6\left((5-1)(2r) + \frac{\pi r}{3}\right)$

**Case-5: $n = 6$ (i.e. six identical circles are tangent to each side of regular hexagon)**

Similarly, consider total $N$ number of identical circles each of radius $r$ the most densely packed in a regular hexagon ABCDEF of side $a$ such that six identical circles are tangent to each side say side ED (as shown in the figure-25(a) below).

In this case, the unknown parameters $N, N_v, a, R,$ and $L$ can be easily determined as given below.

In order to find out the number of identical circles and number of voids, the centres of circles lying on the sides & vertices of internal regular hexagons are counted in CW circular fashion (as shown in the figure-24 below) as follows

$N$ = Total number of identical circles packed in regular hexagon
$$= 1 + 6 + 12 + 18 + 24 + 30$$

It is worth noticing that 6 voids are always present at 6 vertices/corners of the outermost regular hexagon ABCDEF (see figure-24).

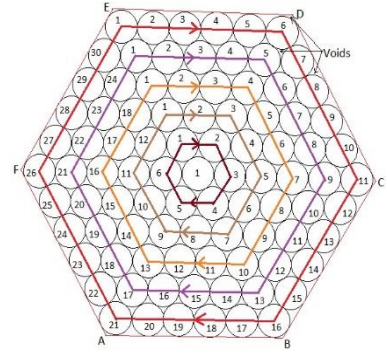

$N_v$ = Total number of voids (unoccupied regions) in regular hexagon
$$= 6 + 2(6 + 12 + 18 + 24 + 30)$$

In right $\Delta C_1 GA$, We have $\angle AC_1G = \frac{\angle AOB}{2} = \frac{60^o}{2} = 30^o$, $C_1G = r$ (see the figure-25(a) below)

$$\tan \angle AC_1G = \frac{AG}{C_1G} \Rightarrow AG = C_1G \tan \angle AC_1G = r \tan 30^o = \frac{r}{\sqrt{3}}$$

**Figure-24: The centers of circles lie on inner hexagons which can be easily counted in CW circular fashion around the central circle.**

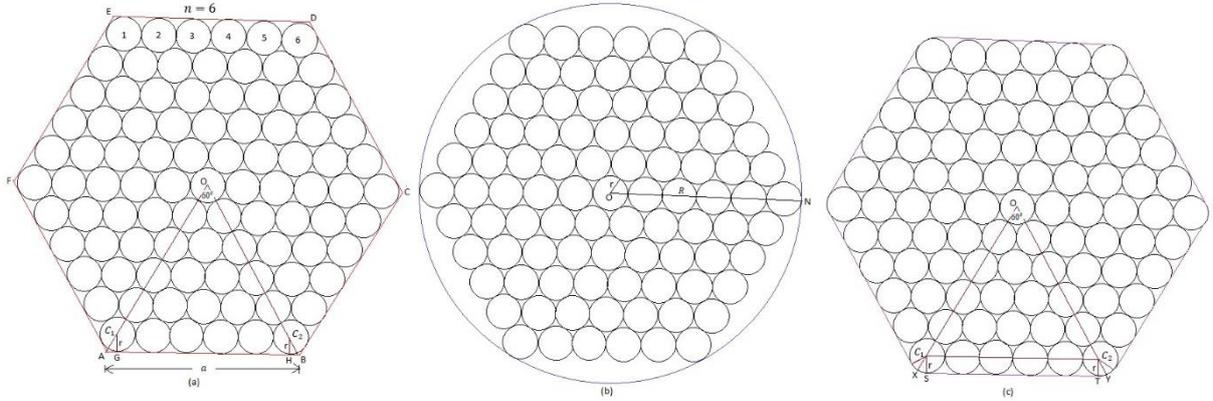

**Figure-25: (a) Six identical circles labeled with 1, 2, 3, 4, 5, 6 are tangent to (each) side ED of regular hexagon ABCDEF (n=6) (b) smallest circle circumscribing all the packed circles (c) minimum length of one turn of flexible & inextensible string wrapped around the identical circles most densely packed in a regular hexagon.**

$$a = \text{Side of regular hexagon ABCDEF} = AB = AG + GH + HB = GH + 2(AG) = 5(2r) + \frac{2r}{\sqrt{3}} = (6-1)(2r) + \frac{2r}{\sqrt{3}}$$

$R$ = Radius of smallest circumscribing circle = $ON = r + 5(2r) = r + (6-1)(2r)$     (see figure 25(b) above)

In the figure-25(c) above, We have $\angle XC_1S = \frac{\angle C_1 O C_2}{2} = \frac{60^o}{2} = \frac{\pi}{6}$, $C_1S = C_1X = r$

$$\Rightarrow \text{Arc } SX = C_1X (\angle XC_1S) = r\left(\frac{\pi}{6}\right) = \frac{\pi r}{6}$$

$$\therefore XY = \text{arc } XS + ST + \text{arc } TY = ST + 2(\text{arc } XS) = 10r + 2\left(\frac{\pi r}{6}\right) = 5(2r) + \frac{\pi r}{3} \quad \text{(see figure 25(c) above)}$$

$$L = \text{Minimum length of one turn of string } = 6(XY) = 6\left(5(2r) + \frac{\pi r}{3}\right) = 6\left((6-1)(2r) + \frac{\pi r}{3}\right)$$

Similarly, consider $n$ number of identical circles touching each side say side AB of regular hexagon ABCDEF then following the symmetry in the above five cases (i.e. $n = 2, 3, 4, 5, 6, ...$) (see the figure-26 below), the unknown parameters $N, N_v, a, R$, and $L$ can be easily determined as follows

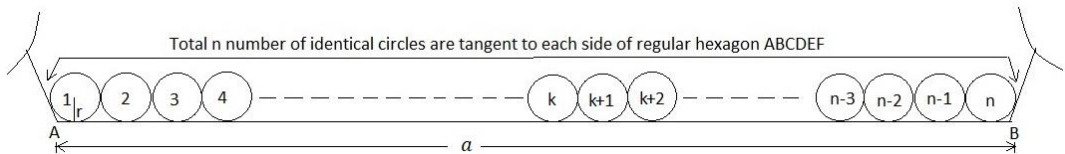

**Figure-26:** $n$ number of identical circles labeled with 1, 2, 3,…k, k+1, …. n are tangent to (each) side AB of regular hexagon ABCDEF (n=n).

Total number of identical circles ($N$) packed in regular hexagon ABCDEF

$$N = 1 + 6 + 12 + 18 + 24 + 30+ \ldots \ldots \text{upto n terms}$$

$$= 1 + 6(1 + 2 + 3 + 4 + 5+ \ldots \ldots \text{upto } (n-1) \text{ terms})$$

$$= 1 + 6\left(\frac{(n-1)(n-1+1)}{2}\right) = 1 + 3n(n-1) = 3n^2 - 3n + 1$$

$$\therefore \ \boldsymbol{N = 3n^2 - 3n + 1} \qquad \ldots \ldots \ldots \ldots (20)$$

Total number of voids (unoccupied regions) ($N_v$) in regular hexagon

$$N_v = 6 + 2(6 + 12 + 18 + 24 + 30+ \ldots \ldots \text{upto } (n-1) \text{ terms})$$

$$= 6 + 12(1 + 2 + 3 + 4 + 5+ \ldots \ldots \text{upto } (n-1) \text{ terms})$$

$$= 6 + 12\left(\frac{(n-1)(n-1+1)}{2}\right) = 6 + 6n(n-1) = 6(n^2 - n + 1)$$

$$\therefore \ \boldsymbol{N_v = 6(n^2 - n + 1)} \qquad \ldots \ldots \ldots \ldots (21)$$

The side ($a$) of the regular hexagon

$$a = (n-1)(2r) + \frac{2r}{\sqrt{3}} = 2r\left(n + \frac{1}{\sqrt{3}} - 1\right)$$

$$\therefore \ \boldsymbol{a = 2r\left(n + \frac{1}{\sqrt{3}} - 1\right)} \qquad \ldots \ldots \ldots \ldots (22)$$

The radius ($R$) of smallest circle circumscribing all the identical circles packed in the regular hexagon

$$R = r + (n-1)(2r) = (2n-1)r$$

$$\therefore \ \boldsymbol{R = (2n-1)r} \qquad \ldots \ldots \ldots \ldots (23)$$

The minimum length ($L$) of the string

$$L = 6\left((n-1)(2r) + \frac{\pi r}{3}\right) = 2(6n + \pi - 6)r$$

$$\therefore \ \boldsymbol{L = 2(6n + \pi - 6)r} \qquad \ldots \ldots \ldots \ldots (24)$$

## 7.1 Packing fraction/density of identical circles most densely packed in a regular hexagon

Consider total $N$ number of identical circles each of radius $r$ most densely packed in a regular hexagon ABCDEF of each side $a$ such that $n$ number of circles are tangent to each side of the hexagon (refer to the figure-25 above).

The total number ($N$) of identical circles most densely packed in regular hexagon is given by the above Eq(20)

$$N = 3n^2 - 3n + 1$$

The side ($a$) of regular hexagon is given by the above Eq(22)

$$a = 2r\left(n + \frac{1}{\sqrt{3}} - 1\right)$$

Now, the packing fraction or packing density ($\rho$) of identical circles in regular hexagon is given as follows

$$\rho = \frac{\text{Area occupied by } N \text{ number of identical circles each of radius } r}{\text{Area of regular hexagon with each side } a}$$

$$= \frac{N(\pi r^2)}{\frac{3\sqrt{3}}{2}a^2} = \frac{(3n^2 - 3n + 1)\pi r^2}{\frac{3\sqrt{3}}{2}\left(2r\left(n + \frac{1}{\sqrt{3}} - 1\right)\right)^2} = \frac{(3n^2 - 3n + 1)\pi}{2\sqrt{3}(\sqrt{3}(n-1) + 1)^2} = \frac{(3n^2 - 3n + 1)\pi\sqrt{3}}{6(\sqrt{3}(n-1) + 1)^2}$$

$$\therefore \rho = \frac{\pi\sqrt{3}}{6}\left(\frac{3n^2 - 3n + 1}{\left(\sqrt{3}(n-1) + 1\right)^2}\right) \quad \ldots\ldots\ldots\ldots\ldots (25)$$

It is clear from the above Eq(25) that the packing fraction/density of identical circles of a finite radius most densely (hexagonally) packed in a regular hexagon depends on the number $n$ of identical circles tangent to each side but independent of the radius of identical circles.

The variation of packing fraction or density ($\rho$) w.r.t. to number of most densely packed identical circles ($n$) tangent to each side of a regular hexagon has been shown in the figure-27.

MATLAB Code to plot variation of packing density ($\rho$)
w.r.t. number of circles ($n$) tangent to each side

```
clc;
clear all;
close all;
syms x y
% x=[0:0.2:6.5];
y=(pi*sqrt(3)/6)*((3*x.^2-3*x+1)/(sqrt(3)*(x-1)+1)^2);
fplot(x,y,[2, 100]);
title('\rho vs n');
xlabel('Number of circles (n) \rightarrow');
ylabel('Packing density (\rho)\rightarrow');
grid on;
```

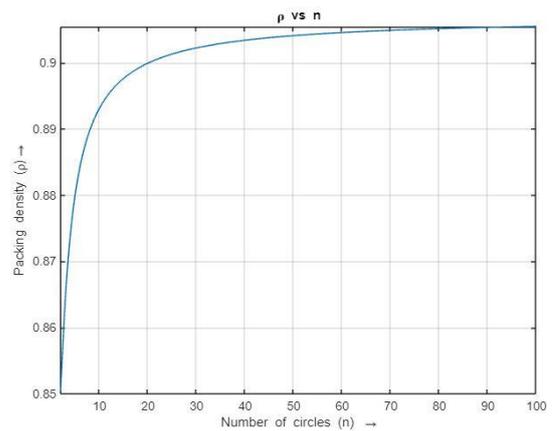

**Figure-27: Variation of packing density ($\rho$) w.r.t. number of most densely packed identical circles ($n$) tangent to each side of regular hexagon.**

It is clear from the above graph that the packing density ($\rho$) of identical circles most densely packed in a regular hexagon increases with the number of identical circles ($n$) tangent to each side and reach to a maximum value when $n \to \infty$.

The minimum value of packing density $\rho_{\min}$ of identical circles most densely packed in a regular hexagon can be determined by substituting $n = 2$ in the above Eq(25) as follows

$$\rho_{\min} = \frac{\pi\sqrt{3}}{6}\left(\frac{3(2)^2 - 3(2) + 1}{\left(\sqrt{3}(2-1)+1\right)^2}\right) = \frac{\pi\sqrt{3}}{6}\left(\frac{7}{\left(\sqrt{3}+1\right)^2}\right) = \frac{\pi\sqrt{3}}{6}\left(\frac{7\left(\sqrt{3}-1\right)^2}{\left(\left(\sqrt{3}+1\right)\left(\sqrt{3}-1\right)\right)^2}\right) = \frac{7\pi\left(2\sqrt{3}-3\right)}{12}$$

$$\therefore \; \rho_{\min} = \frac{7\pi}{12}\left(2\sqrt{3}-3\right) = 0.8505106310376239 \qquad \ldots\ldots\ldots\ldots\ldots (26)$$

$\%\rho_{\min} = 85.05$

It is clear from the above graph that the packing density ($\rho$) has its minimum value $\approx 0.85$ for $n = 2$.

### 7.2 Maximum packing fraction/density of identical circles packed on an infinite plane

In order to obtain the maximum packing fraction or density ($\rho$), the identical circles have to be packed in a regular hexagon. An infinite plane can be assumed to be a regular hexagonal plane with side tending to infinity. In this case, the number of identical circles touching each side of infinite regular hexagon will tend to infinity i.e. $n \to \infty$

Therefore the maximum packing fraction or density of identical circles ($\rho_{\max}$) of a finite radius packed over an infinite plane can be determined by taking limit of packing density $\rho$ given from Eq(25) as $n \to \infty$ as follows

$$\rho_{\max} = \lim_{n\to\infty}\frac{\pi\sqrt{3}}{6}\left(\frac{3n^2-3n+1}{\left(\sqrt{3}(n-1)+1\right)^2}\right) = \lim_{n\to\infty}\frac{\pi\sqrt{3}}{6}\left(\frac{3-\frac{3}{n}+\frac{1}{n^2}}{\left(\sqrt{3}\left(1-\frac{1}{n}\right)+\frac{1}{n}\right)^2}\right) = \frac{\pi\sqrt{3}}{6}\left(\frac{3-0+0}{\left(\sqrt{3}(1-0)+0\right)^2}\right) = \frac{\pi\sqrt{3}}{6}$$

$$\therefore \; \boldsymbol{\rho_{\max} = \frac{\pi\sqrt{3}}{6} = 0.9068996821171088} \qquad \ldots\ldots\ldots\ldots (27)$$

$\%\rho_{\max} = 90.69$

It is clear from the graph in above figure-27 that the packing fraction/density ($\rho$) of most densely packed identical circles in a regular hexagon starts from its minimum value $\approx 0.85$ for $n = 2$ number of identical circles ($n$) tangent to each side of regular hexagon and reaches its maximum value $\rho_{\max} = \frac{\pi\sqrt{3}}{6} \approx 0.9069$ when $n \to \infty$ i.e. an infinite plane.

It is worth noticing that $\boldsymbol{\rho_{\max} = \frac{\pi\sqrt{3}}{6} \approx 0.9069}$ **is maximum possible packing fraction/density of identical circles of a finite radius packed over an infinite plane.**

Alternatively, the maximum packing fraction/density $\rho_{\max}$ of identical circles packed over an infinite plane can be obtained by considering a circle of finite radius $r$ and centre O surrounded by 6 external tangent identical circles with centres A, B, C, D, E and F packed in the hexagonal pattern (as shown in the figure-28). Consider a regular hexagon ABCDEF by joining the centres of 6 identical circles surrounding the identical circle with centre O.

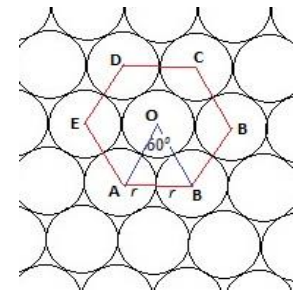

Therefore the maximum packing fraction or density of identical circles ($\rho_{\max}$) is determined as follows

$$\rho_{\max} = \frac{\text{Total area occupied by circles in regular hexagon ABCDEF}}{\text{Area of regular hexagon ABCDEF}}$$

**Figure-28: The identical circles are most densely packed over an infinite plane.**

$$\rho_{\max} = \frac{3(\pi r^2)}{\frac{3\sqrt{3}}{2}(2r)^2} = \frac{\pi}{2\sqrt{3}} = \frac{\pi\sqrt{3}}{6}$$

# 8. 2D packing of circles in the regions bounded by two external tangent circles and the circumscribed circle

When two external tangent circles of radii say $a, \& b$ are inscribed in a third circumscribed circle of radius R then the unoccupied region is divided into two regions which are always unequal (i.e. minor and major) except when $R = a + b$ (See the above figures 2 & 3 or the figure-29 below).

The circle packing in a plane region, bounded by three circles out of which two are external tangent circles touching the third circle internally (as shown in the figure-29 below), can be easily analysed and formulated by substituting $a = r_i$ (i.e. variable/input radius) and $c_{min} = r_{i+1}$ (i.e. radius of next packed circle) in the Eq(1) (as derived in the above section 1.1) to analytically compute the radii of packed circles using the generalized recurrence formula given by as follows

$$r_{i+1} = \frac{r_i bR \left((R-b)r_i + bR - 2\sqrt{bRr_i(R-b-r_i)}\right)}{\left((R+b)r_i + bR\right)^2 - 4r_i bR^2} \qquad \forall \; r_{i+1} < r_i \qquad \ldots\ldots\ldots\ldots(28)$$

Or

$$r_{i+1} = \frac{r_i bR \left((R-b)r_i + bR + 2\sqrt{bRr_i(R-b-r_i)}\right)}{\left((R+b)r_i + bR\right)^2 - 4r_i bR^2} \qquad \forall \; r_i \leq r_{i+1} \leq R-b \qquad \ldots\ldots\ldots\ldots(29)$$

Where, $r_1 = a$ is the radius of one of two external tangent circles i.e. initial circle and $b$ is the radius of second external tangent circle i.e. reference circle. Both the external tangent circles internally touch the third (i.e. largest) circumscribing circle of radius $R$ and divide unoccupied region into two parts i.e. region-1 & region-2 (see figure-29 below).

1) The first recurrence relation i.e. Eq(28) always gives the diminishing or decreasing radii of packed circles in minor region-1 (i.e. labelled with grey colour) and therefore Eq(28) must be used for computing radii of circles packed in minor region-1 for given values of $a, b,$ and $R$ (as shown in the figure-29).

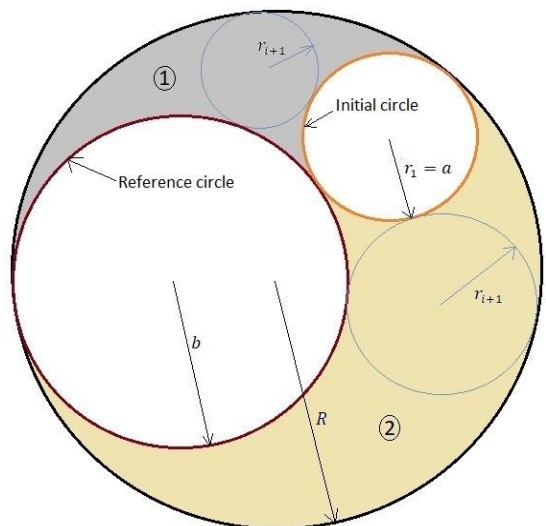

**Figure-29: Circle packing in the plane regions 1 and 2 (i.e. labeled with grey and yellow colors) bounded by three mutually tangent circles with radii $a, b$ and $R$**

2) The second recurrence relation i.e. Eq(29) gives the increasing radii up to a certain maximum value such that $r_{i+1} \leq R - b$ of packed circles in major region-2 (i.e. labelled with yellow colour in the figure-29 above) and after attaining maximum value, the radius $r_{i+1}$ starts oscillating from $r_i$ to its maximum value recursively. Therefore the Eq(29) must be used for packing circles of only increasing radii in major region-2 such that $r_i \leq r_{i+1} \leq R - b$ and after $r_{i+1}$ reaches its maximum value, Eq(28) must be used to further pack required number of circles (of decreasing radii) in the major region-2.

It is worth noticing that in circle packing into major region-2 using the Eq(29),

$$\text{If } r_i = \frac{4bR(R-b)}{(R+b)^2} \Rightarrow r_{i+1} = r_i$$

The above case shows that the radii $r_i$ and $r_{i+1}$ of two circles consecutively packed into major region-2 (i.e. labelled with yellow colour) will be equal and both will be tangent to the diameter line of circumscribing circle of radius $R$ (as shown in the figure-30 below). In this case, the maximum possible value of $r_{i+1}$ is given as

$$(\boldsymbol{r_{i+1}})_{\mathbf{max}} = \frac{\boldsymbol{4bR(R-b)}}{(\boldsymbol{R+b})^2} \qquad \forall\, r_i = \frac{4bR(R-b)}{(R+b)^2}$$

In this case of circles packing into major region-2, after attaining the maximum value i.e. $(r_{i+1})_{max} = \frac{4bR(R-b)}{(R+b)^2}$, the radius $r_{i+1}$ will not oscillate because $r_{i+1} = r_i$ i.e. radii of two consecutively packed circles become equal. This is Non Resonance (NR) condition of the above recursive Eq(29).

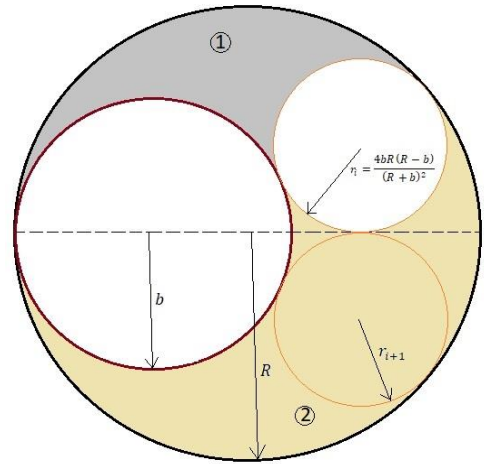

While in all other cases of circle packing in major region-2 using recursive Eq(29), the radius $r_{i+1}$ of packed circle will always oscillate from $r_i$ to $(r_{i+1})_{max} \neq \frac{4bR(R-b)}{(R+b)^2}$ because the radii of two consecutively packed circles are distinct i.e. $r_{i+1} \neq r_i$. This is Resonance (R) condition of the above recursive Eq(29).

i) The value of $r_{i+1}$ obtained from the above recursive Eq(29) will not oscillate if and only if $r_i = \frac{4bR(R-b)}{(R+b)^2}$

**Figure-30:** The radius $r_{i+1}$ of circle packed into major region-2 (i.e. labelled with yellow color) will not oscillate if and only if $r_i = \frac{4bR(R-b)}{(R+b)^2}$

ii) The value of $r_{i+1}$ obtained from the above recursive Eq(29) will always oscillate if $r_i \neq \frac{4bR(R-b)}{(R+b)^2}$

3) In the above recurrence formula i.e. Eq(28) and Eq(29), one of two radii $a$ and $b$ is changed to $r_1$ (while other is kept constant/fixed) and $r_{i+1}$ takes the value of radius of a newly packed circle $r_i$ recursively until the desired number of circles are packed in the bounded plane region-1 or region-2 or both (see figure-29 above).

4) In order to pack $n_1$ number of circles into minor region-1 and $n_2$ number of circles into major region-2 of area unoccupied by two external tangent circles of given radii $a$, and $b$ which internally touch the third circle of radius $R$, the following MATLAB code can be used to determine the radii of all the desired number of circles packed into two unoccupied regions

```
clc;
clear all;
close all;
R=input('Enter the radius of circumscribing circle (in required unit) = ');
a=input('Enter the radius of variant external circle. a = ');
b=input('Enter the radius of reference external circle, b = ');
n1=input('Enter the number of circles to be packed in minor region= ');
n2=input('Enter the number of circles to be packed in major region= ');
```

```matlab
disp('(1) Radii of circles packed into minor region');
x=[];
for i=1:n1
    r(1)=a;
r(i+1)=double((r(i)*b*R*(r(i)*(R-b)+b*R-2*sqrt(r(i)*b*R*(R-b-r(i)))))/((r(i)*(R+b)+b*R)^2-4*r(i)*b*R^2));
A=[x,real(r(i+1))];
disp(A);
end
disp('(2) Radii of circles packed into major region');
for i=1:n2+10
    s(1)=a;
    s(i+1)=double((s(i)*b*R*(s(i)*(R-b)+b*R+2*sqrt(s(i)*b*R*(R-b-s(i)))))/((s(i)*(R+b)+b*R)^2-4*s(i)*b*R^2));
    M=[x,s(i)];
    m1=max(M);

end
for i=1:n2+11
    s(1)=a;
    s(i+1)=double((s(i)*b*R*(s(i)*(R-b)+b*R+2*sqrt(s(i)*b*R*(R-b-s(i)))))/((s(i)*(R+b)+b*R)^2-4*s(i)*b*R^2));
    M=[x,s(i)];
    m2=max(M);
end
m=max(m1,m2);

for i=1:n2
    t(1)=a;
    t(i+1)=double((t(i)*b*R*(t(i)*(R-b)+b*R+2*sqrt(t(i)*b*R*(R-b-t(i)))))/((t(i)*(R+b)+b*R)^2-4*t(i)*b*R^2));

    if t(i+1)>t(i)+10^(-17) && t(i+1)<m+10^(-17)
         n=i;
         B=[x,real(t(i+1))];
         disp(B);
    else
         n=0;
         B=A;
    end
end
    c_max=double((m*b*R*(m*(R-b)+b*R+2*sqrt(m*b*R*(R-b-m))))/((m*(R+b)+b*R)^2-4*m*b*R^2));
    c_min=double((m*b*R*(m*(R-b)+b*R-2*sqrt(m*b*R*(R-b-m))))/((m*(R+b)+b*R)^2-4*m*b*R^2));
disp(m);
disp('Set-1');

for i=1:n2-n
  u(1)=c_max;
  u(i+1)=double((u(i)*b*R*(u(i)*(R-b)+b*R-2*sqrt(u(i)*b*R*(R-b-u(i)))))/((u(i)*(R+b)+b*R)^2-4*u(i)*b*R^2));
  C=[x,real(u(i))];
  disp(C);
end
 disp('Set-2');

for i=1:n2-n

  v(1)=c_min;
  v(i+1)=double((v(i)*b*R*(v(i)*(R-b)+b*R-2*sqrt(v(i)*b*R*(R-b-v(i)))))/((v(i)*(R+b)+b*R)^2-4*v(i)*b*R^2));
  D=[x,real(v(i))];
  disp(D);
end
```

Note: In the above MATLAB code, two sets i.e. Set-1 and Set-2 of radii of circles packed into the major region are obtained out of which only one set of non-repetitive radii is taken into consideration while the other set gives the same values of radii of circles as packed into minor region. Hence care is taken to select correct set of values.

## Conclusions

In this paper, the problems of circles packed into the region bounded by the circular arcs, and the straight lines including square, semi and quarter circles as well as packing of identical circles on infinite plane have been explained and solved in details. The generalized formula have been derived and established for analysing and formulating the circle packing in 3D space as well as on the spherical surface. The simulated values have been compared and verified by analytic formula. The mathematical relations, derived here, can be applied for solving the problems of dense packing of circles and spheres in 2D & 3D containers in industrial applications, and finding the planar density on crystallographic plane.

## Conflict of Interest

I declare no conflicts of interest related to this article.

Harish Chandra Rajpoot

## Data Availability

The data supporting the findings of this study are available and self-sufficient within the article. Raw data that support this study are available from the corresponding author, upon reasonable request.